\documentclass[12pt]{amsart}
\usepackage{a4wide}
\usepackage{eucal}
\usepackage{umlaut}
\usepackage{latexsym}
\usepackage{amssymb}
\usepackage{verbatim}
\usepackage[active]{srcltx}
\usepackage[all]{xy}
\usepackage{listings}
\usepackage[pdfauthor   = {Mohamed Barakat \& Daniel Robertz},
            pdftitle    = {homalg: A meta-package for homological algebra},
            pdfsubject  = {13D02; 13D05; 13D07; 16E05; 16E10; 16E30},
            pdfkeywords = {},
            bookmarks=true,
            bookmarksopen=true,
            colorlinks=true,
            pagebackref=true,
            hyperindex=true,
            linkcolor=blue,
            pagecolor=blue,
            citecolor=blue,
            urlcolor=blue,
            ps2pdf=true
            ]{hyperref}

\usepackage{maple2e}
\LeftMapleSkip=0em

\input amathmoe.sty




\DeclareMathOperator{\Hom}{Hom}
\DeclareMathOperator{\Ext}{Ext}
\DeclareMathOperator{\Tor}{Tor}
\DeclareMathOperator{\coker}{coker}
\DeclareMathOperator{\img}{im}
\DeclareMathOperator{\LL}{L}
\DeclareMathOperator{\RR}{R}

\renewcommand\F{\mathcal{F}}
\newcommand\id{\mathrm{id}}
\renewcommand\phi{\varphi}
\newcommand\CC{\mathcal{C}}
\newcommand\Hull{\mathrm{Hull}}
\newcommand\Obj{\mathtt{Obj}}
\newcommand\Mor{\mathtt{Mor}}

\newcommand\homalg{\mathtt{homalg}}
\newcommand\FunctorMap{\mathtt{FunctorMap}}
\newcommand\Cokernel{\mathtt{Cokernel}}
\newcommand\CokernelMap{\mathtt{CokernelMap}}
\newcommand\Kernel{\mathtt{Kernel}}
\newcommand\KernelMap{\mathtt{KernelMap}}

\newcommand\DefectOfHoms{\mathtt{DefectOfHoms}}
\newcommand\DefectOfHomsMap{\mathtt{DefectOfHomsMap}}
\newcommand\TensorProduct{\mathtt{TensorProduct}}

\newcommand\ttHom{\mathtt{Hom}}
\newcommand\ttHomR{\mathtt{Hom\_R}}
\newcommand\ttHomMap{\mathtt{HomMap}}
\newcommand\ttHomMapR{\mathtt{HomMap\_R}}
\newcommand\ttHomHomR{\mathtt{HomHom\_R}}
\newcommand\ttHomHom{\mathtt{HomHom}}
\newcommand\ttLHomHom{\mathtt{LHomHom}}
\newcommand\ttExt{\mathtt{Ext}}
\newcommand\tA{\mathtt{A}}
\newcommand\tB{\mathtt{B}}
\newcommand\tb{\mathtt{b}}
\newcommand\tC{\mathtt{C}}
\newcommand\tc{\mathtt{c}}
\newcommand\tG{\mathtt{G}}
\newcommand\tL{\mathtt{L}}
\newcommand\tM{\mathtt{M}}
\newcommand\tN{\mathtt{N}}
\newcommand\tS{\mathtt{S}}
\newcommand\tX{\mathtt{X}}
\newcommand\tx{\mathtt{x}}
\newcommand\tY{\mathtt{Y}}
\newcommand\ty{\mathtt{y}}
\newcommand\tz{\mathtt{z}}

\newcommand\cocoa{{\hbox{\rm C\kern-.13em o\kern-.07em C\kern-.13em o\kern-.15em A}}}

\bdoc

\lstset{basicstyle=\tiny,
        frame=single,
        stringstyle=\ttfamily,
        showstringspaces=true}

\author{Mohamed Barakat \& Daniel Robertz}
\address{Lehrstuhl B f\"ur Mathematik, RWTH-Aachen University, 52062 Germany}
\email{\href{mailto:Mohamed Barakat <mohamed.barakat@rwth-aachen.de>}{mohamed.barakat@rwth-aachen.de}, \href{mailto:Daniel Robertz <daniel@momo.math.rwth-aachen.de>}{daniel@momo.math.rwth-aachen.de}}

\PUSH{homalg_package.tex}%
\title{$\homalg$\\ \medskip A meta-package for homological algebra}
\date{2007}

\begin{abstract}
The central notion of this work is that of a functor between categories of finitely presented modules over so-called computable rings, i.e.\  rings $R$ where one can algorithmically solve inhomogeneous linear equations with coefficients in $R$. The paper describes a way allowing one to realize such functors, e.g.\  $\Hom_R$, $\otimes_R$, $\Ext^i_R$, $\Tor^R_i$, as a mathematical object in a computer algebra system. Once this is achieved, one can compose and derive functors  and even iterate this process without the need of any specific knowledge of these functors. These ideas are realized in the ring independent package $\homalg$. It is designed to extend any computer algebra software implementing the arithmetics of a computable ring $R$, as soon as the latter contains algorithms to solve inhomogeneous linear equations with coefficients in $R$. Beside explaining how this suffices, the paper describes the nature of the extensions provided by $\homalg$.
\end{abstract}

\maketitle
%

\section{Introduction}

In the setup of finitely presented module categories, $\homalg$ realizes functors  as mathematical objects which, up to now, can be composed and derived. To this end it realizes, unlike all present systems, the functors not only by the object part, but {\em automatically} also by the morphism part.

$\homalg$ is abstractly designed and therefore can be used as an extension of other mathematical software providing the necessary ring arithmetics in any concrete problem.

\subsection{$\homalg$ as a programming environment for homological algebra}
Functors map objects of a source category to objects of a target category and, in a compatible way, morphisms between two objects in the source category to morphisms between their images in the target category. So when one implements a functor one has not only to take care of how it acts on objects, but also of how it acts on morphisms between these objects.

Homological algebraic constructions \cite{ce,HS,MR,rot,weihom} are present in most of the computer algebra systems such as $\mathsf{Macaulay\  2}$~\cite{M2} , $\mathsf{Singular}$~\cite{singular}/$\mathsf{Plural}$ \cite{plural}, and \cocoa~\cite{cocoa}. One can often find a procedure, let us call it $\ttHom$, to compute $\Hom_R(A,B)$ for two finitely presented modules $A$ and $B$ over commutative rings that are implemented in these systems. Now given two further modules $M$ and $N$ and a morphism $M\xrightarrow{\phi}N$, what about applying $\Hom_R(A,-)$ or $\Hom_R(-,B)$ to the morphism $\phi$? Mathematically it is clear how $\Hom_R(A,-)$ or $\Hom_R(-,B)$ induces a morphism $\Hom_R(A,M)\xrightarrow{\Hom_R(A,\phi)}\Hom_R(A,N)$ or $\Hom_R(N,B)\xrightarrow{\Hom_R(\phi,B)}\Hom_R(M,B)$, but on the level of a computer implementation the information to compute these induced morphisms is normally not contained in the procedure $\ttHom$. So one needs to write or find two completely independent procedures, say $\ttHomMap$ and $\mathtt{Hom2Map}$ implementing the computations of the induced morphisms. $\ttHomMap$ is for example essential\footnote{Maybe not in full generality, since one only needs to apply $\Hom_R(-,B)$ to morphisms between free modules $R^{1\times l_i} \xrightarrow{\phi} R^{1\times l_i'}$.} if one wants to compute $\Ext_R^i(A,B):=\RR^i\Hom_R(-,B)(A)$, the $i^\mathrm{th}$ right derived functor of $\Hom_R(-,B)$ applied to $A$.  Again, a procedure, let us call it $\ttExt$, implementing for all $i$ the computation of the extension modules $\Ext_R^i(A,B)$ is normally not of much help when it comes to computing the induced morphisms $\Ext_R^i(\phi,B)$ or $\Ext_R^i(A,\phi)$. Writing the corresponding procedures, say $\mathtt{ExtMap}$ and $\mathtt{Ext2Map}$, requires one to know about the derivation process, and of course to have $\ttHom$, $\ttHomMap$ and $\mathtt{Hom2Map}$ predefined. So what about composing functors? Consider $\Ext^j_R(\Ext^k_R(-,-),-)$ for example (cf.~\cite{Roos}, \cite[p.~58ff]{Bjo}). For the object part one needs to compose the procedure $\ttExt$ with itself using the first argument. Let us call this compositum $\mathtt{ExtExt}$. The morphism part has now three components. Let us consider $\Ext^j_R(\Ext^k_R(A,-),B)$ for example. Again, a procedure $\mathtt{ExtExt}$ is offhand not of much help. For this morphism part one rather needs a predefined $\mathtt{Ext2Map}$ to first compute $\psi:=\Ext^k_R(A,\phi)$, then a predefined $\mathtt{ExtMap}$ to compute $\Ext^j_R(\psi,B)$. Let us call the resulting procedure $\mathtt{ExtExt2Map}$.

So in order to derive or compose functors, one needs their part on objects but also their part on morphisms. This means, that for each functor one has to implement one procedure for the object part and as many procedures as needed for the morphism part. For more complex functors, constructed out of given ones via iterated compositions and derivations\footnote{Left derivation for covariant functors and right derivation for contravariant functors (=cofunctors).} this quickly becomes unfeasable. So on the level of the computer implementation the following became unavoidable:
\begin{itemize}
  \item Include the mathematical information of how the bifunctor $\Hom_R(-,-)$ acts on
    morphisms inside the procedure $\ttHom$ itself. This has to be done in a way, that a
    general procedure, say $\FunctorMap$, is able to extract this information out
    of $\ttHom$.  Further $\ttHomMap$ and $\mathtt{Hom2Map}$ should now be defined only
    using $\FunctorMap$ applied to $\ttHom$.
  \item Implement a general right derivation procedure for a contravariant functor (given
    by both its parts), let us call it {\tt RightDerivedCofunctor}. Now define
    $\ttExt$ only using {\tt RightDerivedCofunctor} applied to the $\Hom$-functor given by
    both\footnote{Since we right derive $\Hom_R(-,-)$ with respect to its first argument,
    the corresponding morphism procedure $\ttHomMap$ has to be marked in the
    input of {\tt RightDerivedCofunctor}. This is a minor technical issue.} its parts
    [$\ttHom,[\overline{\ttHomMap},\mathtt{Hom2Map}$]]. The procedure $\ttExt$, or rather
    the derivation procedure used to define $\ttExt$, should be able to reconstruct the
    mathematical information of how the bifunctor $\Ext_R(-,-)$ acts on morphisms, and this
    in such a way that the same procedure $\FunctorMap$ mentioned above is able
    to extract this information out of $\ttExt$. The derivation procedure should
    reconstruct this information alone from its input
    [$\ttHom,[\overline{\ttHomMap},\mathtt{Hom2Map}$]]. $\mathtt{ExtMap}$ and
    $\mathtt{Ext2Map}$ should now be defined only using $\FunctorMap$ applied to
    $\ttExt$.
  \item Implement a general composition procedure for two functors (given by both their
    parts), let us call it {\tt Compose\-Functors}. Then define $\mathtt{ExtExt}$ only
    using {\tt Compose\-Functors} applied\footnote{Since we compose $\Ext_R(-,-)$ (here
    with itself) with respect to the first argument, the corresponding morphism
    $\mathtt{ExtMap}$ has to be marked in the first argument of {\tt ComposeFunctors}.} to 
    [$\ttExt,[\overline{\mathtt{ExtMap}},\mathtt{Ext2Map}$]] and
    [$\ttExt,[\mathtt{ExtMap},\mathtt{Ext2Map}$]]. The procedure
    $\mathtt{ExtExt}$, or rather the composition procedure used to define $\mathtt{ExtExt}$,
    should be able to reconstruct the mathematical information of how the trifunctor
    $\Ext^j_R(\Ext^k_R(-,-),-)$ acts on morphisms, in such a way that the same procedure
    $\FunctorMap$ as above is able to extract this information out of
    $\mathtt{ExtExt}$. The composition procedure should reconstruct this information, alone
    out of its input, which is here
    [$\ttExt,[\overline{\mathtt{ExtMap}},\mathtt{Ext2Map}$]] and
    [$\ttExt,[\mathtt{ExtMap},\mathtt{Ext2Map}$]].
    Again, $\mathtt{ExtExtMap}$, $\mathtt{ExtExt2Map}$ and $\mathtt{ExtExt3Map}$ should now
    be defined only using $\FunctorMap$ applied to $\mathtt{ExtExt}$.
\end{itemize}
Like $\FunctorMap$, the procedures {\tt RightDerivedCofunctor} and {\tt Compose\-Functors} should be implemented in a way that is independent of the functors they are applied to. So, with the same {\tt RightDerivedCofunctor} one should be able to define $\RR^i\Ext^j_R(\Ext^k_R(A,-),B)$ and again with the same $\FunctorMap$ compute its part on morphisms, etc.

Starting from Section~\ref{mor} we will try to isolate the mathematical ideas that helped us to realize functors as mathematical objects. That defining more complex functors in {\em both} their parts is now an easy, even automatic task is the {\em first defining property} of $\homalg$.

\subsection{$\homalg$ as a meta-package}\label{meta}
Not a single algorithm to compute in any sort of {\em rings} is implemented in $\homalg$. Rather it is a package built up of homological algebraic definitions and constructions. On the highest levels one finds the construction of connecting homomorphisms and long exact sequences, the processes of composition and derivation of functors and definitions of various specific functors. On lower levels $\homalg$ goes all the way down till it reaches two procedures, which are basically the only ones required from any software implementing the ring specific arithmetics of the (not necessarily commutative) ring $R$. In what follows we refer by ``{\em ring package}'' to such software. To describe the two procedures let $K$ be a (left) $R$-submodule of the free module $R^{1\times n}$ given by finitely many generators. Further let $\tb$ be an arbitrary element of $R^{1\times n}$ and $A$ an $R$-submodule of $K$, again given by finitely many generators ($A\leq K\leq R^{1\times n}\ni \tb$). Let $\tM\in R^{k\times n}$ resp.\  $\tN\in R^{a\times n}$ be the matrix having as rows the generators of $K$ resp.\  of $A$.
\begin{enumerate}
  \item[(Z)] The procedure $\mathtt{DecideZero}(\tb,\tM)$: \\
    Effectively decide if $\tb$ is an element of $K$ or not. ``Effectively'' means:
    In case the element $\tb$ belongs to $K$, then the procedure {\tt DecideZero} returns
    zero and, if asked to, is able to express $\tb$ as an $R$-linear combination of the
    generators of $M$. Otherwise some element $\tb'\in \tb+K$ is returned. This is
    equivalent to deciding the solvability of the inhomogeneous $R$-linear system
    $\tx\tM=\tb$ and, if asked to, finding a particular solution $\tx\in R^{1\times k}$, in
    case one exists. This is simply the straightforward generalization of the ideal
    membership problem to submodules.
  \item[(S)] The procedure $\mathtt{SyzygiesGenerators}(\tM,\tN)$: \\
    Compute a generating set of the $R$-module of solutions of the homogeneous
    $R$-linear system $\tx\tM=0 \mod \tN$. For $\tz\in R^{1\times n}$ the statement
    $\tz=0 \mod \tN$ means that there exists a $\ty\in R^{1\times a}$, such that
    $\tz=\ty\tN$. One calls every such solution a syzygy among the generators of $K$
    modulo $A$.
\end{enumerate}
As the reader may have noticed, it is not required in (Z) that the output $\tb'$ of {\tt DecideZero} only depends on $\tb+K$, i.e.\  $\homalg$ does not require from {\tt DecideZero} to provide a normal form modulo $K$ for every element $\tb\in R^{1\times n}$, except for $\tb\in K$, where {\tt DecideZero} {\em must} return zero. Deciding if two elements $\tb,\tc\in R^{1\times n}$ represent the same class modulo $K$, is reduced to checking if $\tb-\tc\in K$, i.e.\  if $\mathtt{DecideZero}(\tb-\tc,\tM)=0$.

In practice however one often solves (Z) by first constructing a different set of generators for $K$ out of the given one, where this new set satisfies certain properties\footnote{E.g.\   a {\sc Gauss}ian triangular basis in case of fields or an involutive or {\sc Gr\"obner} basis in case of polynomial rings.}, such that a reduction algorithm with respect to a set of generators having these properties is available. In what follows we call such a set of generators a {\em basis} (this should not be confused with a free basis). Normally, such a basis also provides a way to algorithmically solve (S). $\homalg$ enables the ring package to specify a procedure called {\tt BasisOfModule} to compute such a basis. Internally, $\homalg$ only uses the procedure {\tt DecideZero} to perform reductions with respect to such a basis. In case the ring package performs reductions without computing such a basis the procedure {\tt BasisOfModule} is to be set to the identity procedure. In what follows we will refer to the output $\tb'$ of $\mathtt{DecideZero}(\tb,\tM)$ as the {\em reduction} of $\tb$ modulo $K$ (or modulo $\tM$).

Summing up one can say:
\begin{enumerate}
  \item $\homalg$ is designed to be easily extendable by sparing the user the technical
    details of homological constructions.
  \item $\homalg$ is a meta-package that is designed to easily extend mathematical software
    implementing ring arithmetics. All what $\homalg$ needs from such an implementation are
    two procedures: One to effectively solve the ideal membership problem (Z) and one to
    compute a generating set of syzygies, i.e.\  to solve (S) (cf.~\cite{syz,CLO,KR,BTV,DL}). We will call such rings {\em computable}.
\end{enumerate}

We stress the following point: For $\homalg$ it is irrelevant how (Z) and (S) are solved. Hence, the irrelevance of explaining, or even mentioning how to solve these two problems for every particular class of rings is a defining property of $\homalg$, namely its {\em second defining property}.

In Section~\ref{presentation} we describe the categories $\homalg$ is dealing with. In Section~\ref{basic} we deduce from the two procedures {\tt DecideZero} and {\tt SyzygiesGenerators} all the other procedures used in the sequel. After introducing some notation in Section~\ref{cat} we describe $\homalg$'s philosophy of implementing functors in Sections~\ref{mor} and \ref{obj}, whereas Section~\ref{Derivation} describes how $\homalg$ derives functors. Appendix~\ref{RP} outlines the various ring packages that have been successfully used in connection with $\homalg$. Appendix~\ref{implementation} includes some comments to the current implementation. Finally, examples are given in Appendix~\ref{Ex}. All appendices are available at the site of $\homalg$ \cite{homalg}.

It goes without saying that we will suppress some technical issues of the
package for the sake of mathematical clarity. The best (technical) guide to the
package remains its source code, which in nearly all parts can be read as a
mathematical text.

\section{Presentations}
\label{presentation}

Let $R$ be a left noetherian ring with one. Denote by $\CC^\mathrm{mod}$ the category of finitely generated left $R$-modules, which is a full subcategory of the abelian category of left $R$-modules.  And denote by $\CC^\mathrm{mat}$ the category of finite $R$-presentations with objects being finite dimensional matrices over $R$, where one identifies two matrices $\tM\in R^{l_1\times l_0}$  and $\tM'\in R^{r_1\times l_0}$ with the same number $l_0$ of columns to one object, if $R^{1\times l_1}\tM=R^{1\times r_1}\tM'$, as $R$-submodules of $R^{1\times l_0}$.  The set $\mathrm{Mor}_{\CC^{\mathtt{mat}}}(\tM,\tL)$ of morphisms between two objects $\tM\in R^{l_1\times l_0}$ and $\tL\in R^{l_1'\times l_0'}$ is the set  of all $l_0\times l_0'$ $R$-matrices $\phi$ with $R^{1\times l_1}\tM\phi\leq R^{1\times l_1'}\tL$, where one identifies two matrices $\phi_1$ and $\phi_2$ to one morphism, if they induce the same $R$-module homomorphism from $\coker(R^{1\times l_1}\xrightarrow{\tM}R^{1\times l_0})=R^{1\times l_0}/R^{1\times l_1}\tM$ to $\coker(R^{1\times l_1'}\xrightarrow{\tL}R^{1\times l_0'})=R^{1\times l_0'}/R^{1\times l_1'}\tL$.

As usual a presentation is given by generators and relations. If one takes the classes of the $l_0$ standard basis vectors of $R^{1\times l_0}$ as generators of $M:=\coker(R^{1\times l_1}\xrightarrow{\tM}R^{1\times l_0})=R^{1\times l_0}/R^{1\times l_1}\tM$, then the $l_1$ rows of $\tM$ are the defining relations between these generators. We thus call $\tM$ a {\it presentation matrix} or a {\it relation matrix} for $M$. Further we call $R^{1\times l_1}\tM=\img(\tM)$ the {\em relation subspace for $M$ generated by (the rows of)} $\tM$ and denote it by $\langle\tM\rangle$.

If $\tM$ contains a unit in the $j^\mathrm{th}$ column, then the row containing this unit is a relation expressing the $j^\mathrm{th}$ generator of $M=\coker(\tM)$ as an $R$-linear combination in terms of the other generators. Hence one can, using elementary matrix transformations, rewrite the matrix of relations with respect to the remaining generators and obtain a new relation matrix with one less column (and one less row). One iterates this process until the relation matrix is free of units. So, without loss of generality, one can assume that the relation matrix $\tM$ is free of units. For deciding invertability of a ring element and computing its inverse $\homalg$ uses the procedure {\tt Leftinverse} described in \ref{li}.

Summing up, the functor $\coker$ that maps
\[
  R^{l_1\times l_0}\ni\tM\mapsto M:=\coker(R^{1\times l_1}\xrightarrow{\tM}R^{1\times l_0})
\]
is an equivalence between the categories $\CC^\mathrm{mat}$ and $\CC^\mathrm{mod}$. Its inverse functor is to fix a presentation matrix $\tM\in\mathrm{Obj}_{\CC^\mathrm{mat}}$ for each module $M\in\mathrm{Obj}_{\CC^\mathrm{mod}}$ and to use the generators of the fixed presentations to express maps in $\mathrm{Mor}_{\CC^\mathrm{mod}}$ as matrices in $\mathrm{Mor}_{\CC^\mathrm{mat}}$  (cf.~\cite[p.~101]{GP}).

In what follows we no longer distinguish between the two categories and therefore we simply denote both by $\CC$.

\section{Basic procedures}\label{basic}

\subsection{Procedures based on {\tt DecideZero}}

\subsubsection{$\mathtt{RightDivide}(\tB,\tA,\tL)$}\label{rd}
This procedure is the ring-theoretic version of finding a {\it particular solution} of an inhomogeneous linear system of equations. Hence, it is not astonishing that in a lot of computations it plays the decisive role (cf.\  \ref{CompleteImSq}-\ref{pi}).

For three $R$-matrices $\tA,\tB$ and $\tL$ with the same number of {\em columns}, one wants to compute an $\tX$ with
\[
  \tB = \tX\tA \mod \tL,
\]
i.e.\ an $\tX$, such that there exists a $\tY$ with
\[
  \tB=\tX\tA+\tY\tL
  =(\tX,\tY)\left(\begin{array}{c}\tA\\ \tL\end{array}\right).
\]
After computing $(\tX,\tY)$ one often simply throws away $\tY$. So without loss of  generality we forget about $\tL$, i.e.\ consider the system $\tB=\tX\tA$. To find $\tX$ one computes a basis $\tG$ for $\tA$ (i.e.\ a matrix $\tG$ with rows being a basis for the rows of $\tA$) together with a matrix $\tC_\tA$, such that $\tG=\tC_\tA \tA$. Then one reduces $\tB$ modulo $\tG$, i.e.\ computes a matrix $\tN:=\mathtt{DecideZero}(\tB,\tG)$ (cf.~\ref{meta},(Z)), with rows being those of $\tB$ reduced modulo $\tG$. We compute $\tN$ together with a matrix $\tC_\tB$, such that $\tN=\tB+\tC_\tB \tG$. If $\tN$ is not the zero matrix, then the system is not solvable. Otherwise $\tX=-\tC_\tB \tC_\tA$ is a solution. Thinking of $\tX$ as ``$\tX=\tB\tA^{-1}$'' justifies the name of the procedure.

The most prominent application of {\tt RightDivide} in $\homalg$ is the following situation: \\
Given the three modules and the two morphisms
\[
  \xymatrix{
    M' \ar[rd]^\gamma \\
    N'\ar[r]_\beta & N
  }
\]
satisfying the so-called {\it image condition} $\img\gamma\leq\img\beta$ (as submodules of $N$),
find a third morphism $\psi:M'\to N'$ that makes the triangle commute. The image condition is obviously a necessary condition for such a $\psi$ to exist. There are two instances of special importance for us in which such a $\psi$ always exists:
\begin{enumerate}
  \item  $M'$ is a free module (of finite rank $r$), or
  \item $\beta$ is injective.
\end{enumerate}
In the first case let $(b_1,\ldots,b_r)$ be a free basis of $M'$. Define $b_i\psi$ to be any element of the set $(b_i\gamma)\beta^{-1}$, i.e.\ any element in the preimage of $b_i\gamma$ under $\beta$, and since $M'$ is free on $(b_1,\ldots,b_n)$ this extends by linearity to a morphism $\psi:M'\to N'$ satisfying $\psi\beta=\gamma$. In the second case $\beta$ is an isomorphism onto $\img\beta$ and one defines $\psi:=\gamma\beta^{-1}$, where $\gamma$ is viewed as a morphism $M'\to\img\beta$. \\
Now let $\tM',\tN'$ and $\tN$ be relation matrices of $M',N'$ and $N$ respectively. Then we are looking for a matrix $\psi$, together with matrices $Y$ and $Z$ such that
\begin{enumerate}
  \item[(i)] $\gamma=\psi\beta+Y\tN$,
  \item[(ii)] $\tM'\psi=Z\tN'$.
\end{enumerate}
If one regards $\psi$ merely as a matrix, i.e.\  not necessarily a morphism $M'\to N'$, then condition (i) is nothing else but the above mentioned image condition. Therefore, for the situation described above, such a matrix $\gamma$ always exists. But now condition (ii) states that $\psi$ carries the relations of $M'$ to relations of $N'$, i.e.\  it requires $\psi$ to be a morphism from $M'$ to $N'$. While condition (i) is already in the form needed by ${\tt RightDivide}$\footnote{Following the notation used in defining {\tt RightDivide} we set $\tA=\beta$, $\tB=\gamma$ and $\tL=\tN$.}, $\psi$ on the left hand side of condition (ii) is multiplied from the right\footnote{If the ring $R$ is commutative, one can use the {\sc Kronecker} product to construct for any matrix $C$ over $R$ a matrix $\widetilde{C}$ satisfying $\mathrm{row}(CX)=\mathrm{row}(X)\widetilde{C}$ for any composable matrix $X$, where $\mathrm{row}(C)$ is the row vector consisting of the rows of $C$ written behind each other in the obvious order. By this trick one can rewrite the two conditions (i) and (ii) in a single affine condition, where $\mathrm{row}(\psi), \mathrm{row}(Y)$ and $\mathrm{row}(Z)$ are multiplied from the left as required by $\mathtt{RightDivide}$. Since the resulting affine system is, in general, much bigger compared with the initial ones solving it is computationally expensive. Cf.~\cite{zl}.}. It now happens that in the two instances (1) and (2) condition (ii) is automatically fulfilled, whenever (i) is. In the first instance one needs to additionally require that $M'$ is given on a free basis. Then the matrix of relations $\tM'$ is the zero matrix (with one row and $r$ columns), so (ii) is trivially satisfied for any $\psi$ by taking $Z$ to be the zero matrix of appropriate dimension. Case (2) is equivalent to saying that the preimage of the subspace $\langle\tN\rangle$ of relations for $N$ (generated by $\tN$) under $\beta$ coincides with the subspace $\langle\tN'\rangle$ of relations for $N'$ (generated by $\tN'$): $(\langle\tN\rangle)\beta^{-1}=\langle\tN'\rangle$. But since $\gamma$ is a morphism $M'\to N$ any relation $m'\in\langle\tM'\rangle$ of $M$ satisfies $0=m'\gamma\stackrel{\text{(i)}}{=}(m'\psi)\beta \mod \langle\tN\rangle$. Hence, there exists an $n'\in\langle\tN'\rangle$, such that $n'=m'\psi$. So $\psi$ takes relations of $M'$ to relations of $N'$ and is hence a morphism $M'\to N'$.

\subsubsection{$\mathtt{CompleteImSq}$}\label{CompleteImSq}
This is the most prominent reincarnation of {\tt Right\-Divide} in $\homalg$. We call an incomplete square of the form
\[
  \xymatrix{
    M' \ar[r]^\alpha & M\ar[d]^{\phi} \\
    N'\ar[r]^\beta & N
  }
\]
an {\em image square} if the image condition $\img(\alpha\phi)\leq\img(\beta)$ (as submodules of $N$) is satisfied. This is precisely the above situation for $\gamma=\alpha\phi$. In Section~\ref{mor} {\tt CompleteImSq} is applied to image squares with injective $\beta$ (case (2)) and in Subsection~\ref{rm} to ones where $M'$ is free (case (1)). In these two instances, as shown above, the image square is completable by a morphism $\psi:M'\to N'$ which is directly computable using {\tt RightDivide}.

\subsubsection{$\mathtt{Leftinverse}$}\label{li}
The typical applications of {\tt Leftinverse} correspond to the two situations (1) and (2) described for {\tt RightDivide} in \ref{rd} (cf.~Example \ref{ExOre}):
\begin{enumerate}
  \item Either take $M'=N$ {\em free} (case~\ref{rd},(1)), $\gamma=\mathrm{id}$ and
    $\beta:N'\to N$ surjective. Then the image condition is trivially satisfied.
    $\psi:M'=N\to N'$ is then nothing else but the left inverse or a {\em split} of $\beta$.
  \item Or take $M'=N$, $\gamma=\mathrm{id}$ and $\beta:N'\to N$ an isomorphism. Then
    in particular $\beta$ is injective (case \ref{rd},(2)), and the image condition is
    trivially satisfied. $\psi:M'=N\to N'$ is then nothing else but the (left)inverse of
    $\beta$.
\end{enumerate}
 
\subsubsection{$\mathtt{Preimage}$}\label{pi}
Let $M\xrightarrow{\tA}L$ be a morphism between the two modules $M$ and $L:=\coker(\tL)$. If the matrix $\tB$ consists of rows that are in the image of $\tA$, i.e.\  if the image condition is satisfied, then the rows of the matrix $\tX$ computed by {\tt RightDivide} are the preimages of the rows of $\tB$ under $\tA$. In Subsection \ref{ch} we use the procedure {\tt Preimage} to construct connecting homomorphisms.

\subsection{Procedures additionally based on {\tt SyzygiesGenerators}}\label{sg}

For a matrix $\tA\in R^{l_1\times l_0}$ a {\em matrix of generating syzygies} or a {\em syzygy matrix} is a matrix $\tX\in R^{l_2\times l_1}$ such that $R^{1\times l_2}\tX=\ker(R^{1\times l_1}\xrightarrow{\tA} R^{1\times l_0})$. A slight generalization of this is when one requires that $R^{1\times l_2}\tX=\ker(R^{1\times l_1}\xrightarrow{\tA} \coker(\tL))$, where $\tL\in R^{a\times l_0}$ has the same number $l_0$ of columns as $\tA$. One then calls $\tX$ a {\em syzygy matrix} of $\tA$ {\em modulo} $L$, or {\em relative} to the module $\coker(L)$.

So the rows of the matrix $\tX$ of generating syzygies generate the solution space of the homogeneous linear system $\mathtt{x}\tA=0 \mod \tL$. For $\tz\in R^{1\times l_0}$ the statement $\tz=0 \mod \tL$ means that there exists a $\ty\in R^{1\times a}$, such that $\tz=\ty\tL$.

\subsubsection{$\mathtt{ResolutionOfModule}(\tM)$}\label{res}
By iterating the process of taking syzygies one obtains a free resolution of the module $M=\coker(\tM)$:
\[ \cdots\xrightarrow{\phi_{i+2}}
   R^{1\times l_{i+1}}\xrightarrow{\phi_{i+1}} R^{1\times l_{i}}\xrightarrow{\phi_{i}}
   R^{1\times l_{i-1}}\xrightarrow{\phi_{i-1}}
   \cdots \xrightarrow{\phi_{2}} R^{1\times l_{1}}\xrightarrow{\phi_{1}=\tM}
   R^{1\times l_{0}}\to M\to 0
\]
of desired length. What $\homalg$ additionally does is the following: As explained in Section~\ref{presentation} we can assume the relation matrix (or the first syzygy matrix) $\tM=\phi_1$ free of units. Starting from $i=2$, whenever the $i^\mathrm{th}$ syzygy matrix $\phi_i$ is computed, $\homalg$ uses the units appearing in it to locate the redundant rows of the $(i-1)^\mathrm{st}$ syzygy matrix $\phi_{i-1}$. This means, if $\phi_i$ contains a unit in the $j^\mathrm{th}$ column, then the row in $\phi_i$ containing this unit says that the $j^\mathrm{th}$ row of $\phi_{i-1}$ is an $R$-linear combination of the other rows (of $\phi_{i-1}$), so it is redundant for generating the $(i-1)^\mathrm{st}$ syzygies and can be omitted from $\phi_{i-1}$. The $j^\mathrm{th}$ row is thus deleted from $\phi_{i-1}$ and the $i^\mathrm{th}$ syzygy matrix $\phi_i$ is recomputed. The rows of the latter are again used to locate further redundant rows of $\phi_{i-1}$. This process obviously stabilizes. Then one proceeds with the $(i+1)^\mathrm{st}$ syzygy matrix $\phi_{i+1}$, etc. We end up having a free resolution where all the matrices $\phi_i$ are free of units.

In the case of graded modules over a positively graded ring with degree 0 part a field, or modules over local rings, this process indeed yields a minimal free resolution (cf.~\cite[p.~472]{eis} and \cite[p.~40]{Sch}). In both cases eliminating units from the matrices $\phi_i$ of the resolution is equivalent to requiring that all entries belong to the unique maximal graded resp.\  unique maximal ideal of $R$. 

However, and although in general this process does not yield a ``minimal''\footnote{Such a notion is in general not even well defined, cf.~\cite[p.~472]{eis}.} free resolution, it often enough reduces the involved dimensions of the matrices considerably (cf.~\cite{KR2}).

\subsubsection{$\mathtt{SubfactorModule}(\tM_1,\tM_2)$}\label{subfac}
Given two matrices $\tM_1$ and $\tM_2$ with the same number of columns, computing a presentation matrix for the  subfactor module $(\langle\tM_1\rangle+\langle\tM_2\rangle)/\langle\tM_2\rangle$ (short-hand: $\langle\tM_1\rangle/\langle\tM_2\rangle$) goes as follows: One computes a basis $\tB$ of $\tM_2$. Then one reduces $\tM_1$ modulo $\tB$ and gets $\tN$. The syzygy matrix $\tS$ of $\tM$ modulo $\tN$ is the desired presentation matrix.

\section{Categories of complexes of given finite length}\label{cat}

As always, let $\CC$ be the module category of Section~\ref{presentation}. By $D^k(\CC)$ we denote the category of chain complexes $C:C_k\to C_{k-1}\to\cdots\to C_1$ of $\CC$ of finite length $k-1$ and their chain maps.

\bdf[The category $D^k_i(\CC)$]
  For $i\in\{1,\ldots,k\}$ let $D^k_i(\CC)$ be the factor category of
  $D^k(\CC)$ defined by forgetting in a chain map the morphisms between
  all but the $i^\mathrm{th}$ chain module.
\edf
Note that for a morphism in $D^k_i(\CC)$
\[
  \xymatrix@C=1.5pc@R=1.5pc{
    C_k' \ar[r] & \cdots \ar[r] & C_{i+1}' \ar[r] & C_i' \ar[r] \ar[d]^\phi &
    C_{i-1}' \ar[r] & \cdots \ar[r] & C_1' \\
    C_k \ar[r] & \cdots \ar[r] & C_{i+1} \ar[r] & C_i \ar[r] &
    C_{i-1} \ar[r] & \cdots \ar[r] & C_1
  }
\]
there exists at least one commutative completion
\[
  \xymatrix@C=1.5pc@R=1.5pc{
    C_k' \ar[r] \ar[d]^{\phi_k} & \cdots \ar[r] & C_{i+1}' \ar[r] \ar[d]^{\phi_{i+1}} &
    C_i' \ar[r] \ar[d]^\phi & C_{i-1}' \ar[r] \ar[d]^{\phi_{i-1}} &
    \cdots \ar[r] & C_1' \ar[d]^{\phi_1}\\
    C_k \ar[r] & \cdots \ar[r] & C_{i+1} \ar[r] & C_i \ar[r] &
    C_{i-1} \ar[r] & \cdots \ar[r] & C_1
  }
\]
which is a preimage in $D^k(\CC)$.

By definition $D^1_1(\CC)$ is simply a different notation for $\CC$. From now on all functors we consider take values in the module category $\CC$.

\bdf[The functor $\Obj^k_i$]\label{Obj}
  Denote by $\Obj^k_i:D^k_i(\CC)\to \CC$ the full functor mapping a complex to its $i^\mathrm{th}$ module:
  \[
    \xymatrix@C=1.5pc@R=1.5pc{
      C_k' \ar[r] & \cdots \ar[r] & C_{i+1}' \ar[r] & C_i' \ar[r] \ar[d]^\phi &
      C_{i-1}' \ar[r] & \cdots \ar[r] & C_1' \ar@{}[rrd]|{\mapsto} && C_i' \ar[d]^\phi \\
      C_k \ar[r] & \cdots \ar[r] & C_{i+1} \ar[r] & C_i \ar[r] &
      C_{i-1} \ar[r] & \cdots \ar[r] & C_1 && C_i
    }
  \]
\edf

\bdf[The functor $\Mor^k_i$]\label{Mor}
  Denote by $\Mor^k_i:D^k_i(\CC)\to \CC$ the full functor mapping a complex to its
  $(i+1)^\mathrm{st}$ morphism:
  \[
    \xymatrix@C=1.5pc@R=1.5pc{
      C_k' \ar[r] & \cdots \ar[r] & C_{i+1}' \ar[r]^{\mu_{i+1}'} &
      C_i' \ar[r] \ar[d]^\phi & C_{i-1}' \ar[r] & \cdots \ar[r] &
      C_1' \ar@{}[rrd]|{\mapsto} && C_{i+1}' \ar[r]^{\mu_{i+1}'} & C_i'
      \ar[d]^\phi \\
      C_k \ar[r] & \cdots \ar[r] & C_{i+1} \ar[r]^{\mu_{i+1}} & C_i \ar[r] &
      C_{i-1} \ar[r] & \cdots \ar[r] & C_1 && C_{i+1}\ar[r]^{\mu_{i+1}} & C_i
    }
  \]
\edf

\section{The morphism part of a functor}\label{mor}
A functor is by definition a map between two categories that maps objects to objects and, in a compatible way, morphisms to morphisms. To be allowed to speak about functors one needs to be able, not only to define the functor on objects, but also on morphisms between objects.

It might be puzzling for the reader that we start with a section devoted to the morphism part rather than the object part of functors. The reason for that will become clear towards the end of this section.

Here we write everything for covariant functors. Adapting things for contravariant
functors is done in the obvious way.

For a morphism $S \xrightarrow{\phi} T$ and a functor $F$ we want to compute $F(\phi)$. There are three cases $\homalg$ distinguishes:
\begin{enumerate}
  \item[(Cmp)] The functor $F$ is defined in $\homalg$ as a composition of two functors:
    $F=F_1\circ F_2$. 
  \item[(Der)] The functor $F$ is defined in $\homalg$ as the $i^\mathrm{th}$ (left)
    derivation\footnote{$\homalg$ provides procedures to left derive covariant
    functors and right derive contravariant functors. These are the cases computed
    via a projective resolution of the module.} of another functor:
    $F=\LL_i G$.
  \item[(Bsc)] The functor $F$ is defined in  $\homalg$ neither by composition nor by
    derivation. In what follows such functors are called {\em $\homalg$-basic}\footnote{The
    prefix $\homalg$ indicates that the notion of ``basic functor'' is {\em not} a
    mathematical definition. If, for example, the functor $\Hom_R(\Hom_R(-,R),R)$ would
    have been implemented without using $\homalg$'s composition
    procedure {\tt ComposeFunctors} (see \ref{ComposeFunctors}), we then would call it $\homalg$-basic.} functors.
\end{enumerate}
Dealing with (Cmp), i.e.\  with a composed functor is easy: Since $F(\phi)=F_1(F_2(\phi))$ one is able to reduce computing $F(\phi)$ to computing $F_2(\phi)$ and then $F_1(F_2(\phi))$. A bit more involved is the case (Der), i.e.\  when $F=\LL_iG$. In Subsection~\ref{rm} it is shown how to reduce the computation of $F(\phi)$ to essentially computing $G(\phi)$. For $\homalg$-basic functors the idea is to reduce the computation of $F(\phi)$ to completing an image square (cf.~\ref{CompleteImSq}, case (2)). To this end one {\em embeds} $F(S)$ resp.\  $F(T)$ in a module $\Hull_F(S)$ resp.\  $\Hull_F(T)$ with an induced morphism $\Hull_F(\phi)$ determined by $\phi$:
\begin{equation}\tag{Hull}\label{Hull}
  \xymatrix{
    0 \ar[r] & F(S) \ar[d]_{F(\phi)} \ar[r]^<(.2){\iota_F(S)} &
    \Hull_F(S) \ar[d]^{\Hull_F(\phi)} \\
    0 \ar[r] & F(T) \ar[r]^<(.2){\iota_F(T)} & \Hull_F(T)
  }
\end{equation}
Expressed categorically, $\Hull_F$ is a functor, which we call a {\em hull functor} of $F$, and $\iota_F$ is a natural transformation, which we call the {\em corresponding natural embedding}.

For a given functor $F$ the idea is to either define the hull functor $\Hull_F$ from scratch (e.g.\  for\footnote{Since one can directly provide the action of $\Cokernel$ and $\TensorProduct$ on morphisms, one takes them as their own hull functors.} $\Cokernel$ (in \ref{coker}), $\ttHomR$ (in \ref{HomR}), $\TensorProduct$ (in \ref{tp}) and $\ttHom$ (in \ref{hom})) or to define $\Hull_F$ by composition of already defined functors and the forgetful functors $\Obj^k_i$ (Def.~\ref{Obj}) and $\Mor^l_j$ (Def.~\ref{Mor}) (e.g.\  for $\Kernel$ (in \ref{ker}) and $\DefectOfHoms$ (in \ref{def})).

The procedure in $\homalg$ that precisely accomplishes the above mentioned distinction is called $\FunctorMap$ (cf.~Appendix \ref{FunctorMap}).

\section{Some $\homalg$-basic functors}\label{obj}
Since we need to explain how to construct for each functor its hull functor, we need to first of all mention which specific standard method for computing the object part we use (see for example \cite{GP}) and then how to construct the hull functor in the specific setup. The hull functor is the missing piece of data that allows one to automatically compute the functor 
part on morphisms.

All functors we consider are of the form $F:D^k_i(\CC)\to\CC$. We simply distinguish two types of functors, depending on whether $F(C)$ is a subfactor module of $\Obj^k_i(C)$, for all objects $C\in D^k_i(\CC)$, or not. In what follows a functor of the previous kind is called a {\em functor producing subfactor modules}.

To avoid confusing the two parts of a functor with source category\footnote{This is the case in \ref{coker}, \ref{ker}.} $D^2_k(\CC)$ (which has as its set of objects also morphisms between modules), we use two different names for the two parts: If the object part of a functor is called {\tt F}, then its morphism part is called {\tt FMap}. The name {\tt F} is also used to refer to the functor in both its parts.

\subsection{Functors producing subfactor modules}\label{fpsfm}

\subsubsection{The functor $\Cokernel$}\label{coker}
On the level of objects the covariant cokernel functor associates to a morphism between two modules its cokernel. $\Cokernel:D^2_1(\CC)\to\CC$, i.e.\
\[
  \xymatrix@R=1.5pc{
    A' \ar[r]^\alpha & A \ar[d]_\phi \ar@{}[drr]|{\mapsto} &&
    \Cokernel(\alpha) \ar[d]^{\CokernelMap(\phi)} \\
    B' \ar[r]^\beta & B &&
    \Cokernel(\beta)
  }
\]
such that
\[
  \xymatrix@R=1.5pc{
    A' \ar[r]^\alpha & A \ar[r] \ar[d]_\phi &
    \Cokernel(\alpha) \ar[d]^{\CokernelMap(\phi)} \ar[r] & 0\\
    B' \ar[r]^\beta & B \ar[r] &
    \Cokernel(\beta) \ar[r] & 0
  }
\]
is commutative and exact.

Defining the object part $\Cokernel(\alpha)$ is simple: After fixing a presentation matrix $\tA$ for $A$ and generators for $A'$ one can view $\alpha$ as a matrix with the same number of columns as $\tA$. Now take the supermatrix $\left(\begin{array}{c} \alpha \\ \tA \end{array}\right)$ as a presentation matrix for $\Cokernel(\alpha)$.

The cokernel functor is the most basic functor in our setting, in the sense that computing its morphism part $\CokernelMap(\phi)$ is trivial: If we take the residue classes of the generators of $A$ resp.\  $B$ to be the generators of $\Cokernel(\alpha)$ resp.\  $\Cokernel(\beta)$, then the matrix representing $\CokernelMap(\phi)$ coincides with that of $\phi$. Hence one can set the hull functor $\Hull_\Cokernel=\Cokernel$ and the natural embedding $\iota_\Cokernel$ is the identity transformation\footnote{This is not the whole truth. The functor $\Cokernel$ calls a procedure named $\mathtt{Presentation}$, which, among other things, tries with the help of the procedure $\mathtt{BetterGenerators}$, to reduce the number of generators of the resulting module using either normal form algorithms referred to in the introduction and at the end of Appendix~\ref{RP}, in case they are available and applicable, or otherwise, beside the one described in Section~\ref{presentation}, several clever heuristics. Taking care of this, is a simple, but technical issue.}.

\subsubsection{The functor $\Kernel$}\label{ker}
On the level of objects the covariant kernel functor associates to a morphism between two modules its kernel. $\Kernel:D^2_2(\CC)\to\CC$, i.e.\
\[
  \xymatrix@R=1.5pc{
    A \ar[r]^\alpha \ar[d]^\phi & A'' \ar@{}[drr]|{\mapsto} &&
    \Kernel(\alpha)\ar[d]^{\KernelMap(\phi)} \\
    B \ar[r]^\beta & B'' &&
    \Kernel(\beta)
  }
\]
such that
\[
  \xymatrix@R=1.5pc{
    0 \ar[r] & \Kernel(\alpha) \ar[d]_{\KernelMap(\phi)} \ar[r] &
    A \ar[r]^\alpha \ar[d]^\phi & A''\\
    0 \ar[r] & \Kernel(\beta) \ar[r] &
    B \ar[r]^\beta & B''
  }
\]
is commutative and exact.

Defining the object part $\Kernel(\alpha)$ goes like this: After fixing a presentation matrix $\tA$ resp.\  $\tA''$ for $A$ resp.\  $A''$ one can view $\alpha$ as a matrix with the same number of columns as $\tA''$. Compute the syzygy matrix $\mathtt{iota}$ of $\alpha$ modulo $\tA''$ (cf.~the beginning of Subsection \ref{sg}). Now use {\tt SubfactorModule} to compute $\Kernel(\alpha):=\langle\mathtt{iota}\rangle/\langle\tA\rangle$ (cf.~\ref{subfac}). See also \cite[p.~97]{GP}.

In the case of the functor $\Kernel$ the hull functor is obviously $\Hull_\Kernel=\Obj^2_2$ and the natural embedding $\iota_\Kernel$ is the transformation embedding the kernel of a morphism into its source module:
\[
  \xymatrix@1{
   0 \ar[r] & \Kernel(\alpha) \ar[rr]^<(.24){\iota_{\Kernel}(\alpha)} &&
   A=\Obj^2_2(A\stackrel{\alpha}{\to} A'')
  }
\]
The matrix of the natural embedding is simply $\mathtt{iota}$.

\subsubsection{The functor $\DefectOfHoms$}\label{def}
On the level of objects the covariant defect functor associates to two composable
morphisms $\alpha_1$ and $\alpha_2$ with $\alpha_2\alpha_1=0$ their defect of
exactness $\ker(\alpha_2)/\img(\alpha_1)$. $\DefectOfHoms:D^3_2(\CC)\to\CC$, i.e.\
\[
  \xymatrix@R=1.5pc{
    A' \ar[r]^{\alpha_1} & A \ar[r]^{\alpha_2} \ar[d]^\phi &
    A'' \ar@{}[drr]|{\mapsto} & &
    \DefectOfHoms(\alpha_1,\alpha_2)\ar[d]^{\DefectOfHomsMap(\phi)} \\
    B' \ar[r]^{\beta_1} & B \ar[r]^{\beta_2} & B'' &&
    \DefectOfHoms(\beta_1,\beta_2)
  }
\]
such that
\[
  \xymatrix@R=1.5pc{
    0 \ar[r] & \DefectOfHoms(\alpha_1,\alpha_2) \ar[d]_{\DefectOfHomsMap(\phi)}\ar[r]
    & \Cokernel(\alpha_1) \ar[r]^<(.3){\alpha_2} \ar[d]^{\CokernelMap(\phi)} & A''\\
    0 \ar[r] & \DefectOfHoms(\beta_1,\beta_2) \ar[r] &
    \Cokernel(\beta_1) \ar[r]^<(.3){\beta_2} & B''
  }
\]
is commutative and exact.

The definition of the object part of $\DefectOfHoms$ uses the ideas in \ref{coker} and \ref{ker}: Fix a presentation matrix $\tA$ for $A$ and $\tA''$ for $A''$. Fix generators for $A'$. Then one can view $\alpha_1$ resp.\  $\alpha_2$ as a matrix with the same number of columns as $\tA$ resp.\  $\tA''$. Using $\tA$, $\alpha_2$ and $\tA''$ compute the matrix $\mathtt{iota}$ of the embedding of $\Kernel(\alpha_2)$ in $A$ as in \ref{ker}. Then use {\tt SubfactorModule} to compute the defect $\DefectOfHoms(\alpha_1,\alpha_2):=\langle\mathtt{iota}\rangle/\left\langle\left(\begin{array}{c}\alpha_1 \\ \tA \end{array}\right)\right\rangle$.

In the case of the functor $\DefectOfHoms$ the hull functor is obviously $\Hull_\DefectOfHoms=\Cokernel\circ\Mor^3_2$ and the natural embedding $\iota_\DefectOfHoms$ is the transformation embedding the defect of two composable morphisms into the cokernel of the first morphism:
\[
  \xymatrix@1{
    0 \ar[r] & \DefectOfHoms(\alpha_1,\alpha_2)
    \ar[rrr]^<(.31){\iota_{\DefectOfHoms}(\alpha_1,\alpha_2)} &&&
    \Cokernel(\alpha_1),
  }
\]
where $\Cokernel(\alpha_1)=    \Cokernel(\Mor^3_2(A'\stackrel{\alpha_1}{\to}A\stackrel{\alpha_2}{\to} A''))$. The matrix of the natural embedding is $\mathtt{iota}$.

The functors $\Cokernel$ and $\Kernel$ are obviously special cases of this functor.

\subsection{Other types of $\homalg$-basic functors}
Contrary to functors producing subfactor modules, where the morphism part of the hull functor is given by the induced morphism on the factor module, here we need to explicitly mention how the morphism part of the hull functor is defined.

For the rest of the subsection let
\begin{eqnarray*}
  M=\coker(\tM)=R^{1\times l_0}/R^{1\times l_1}\tM,\quad
  L=\coker(\tL)=R^{1\times l_0'}/R^{1\times l_1'}\tL,
\end{eqnarray*}
be two finitely presented modules:
\begin{eqnarray*}
  \xymatrix@1{
    R^{1\times l_1} \ar[r]^{\tM} & R^{1\times l_0} \ar[r]^{\nu_M} & M \ar[r] & 0,
  }\quad
  \xymatrix@1{
    R^{1\times l_1'} \ar[r]^{\tL} & R^{1\times l_0'} \ar[r]^{\nu_L} & L \ar[r] & 0.
  }
\end{eqnarray*} 

Further let $M\stackrel{\phi}{\to}N$ be a morphism, between two finitely presented modules $M$ and
\[
  N=\coker(\tN)=R^{1\times l_0''}/R^{1\times l_1''}\tN,
\]
i.e.\
\[
  \phi\in R^{l_0\times l_0''}.
\]
Note that we write the morphisms on the right, i.e.\ we use the row convention.

\subsubsection{The functor $\ttHomR$}\label{HomR}
Let $R$ be a (not necessarily commutative) ring with a fixed {\em involution}\footnote{An order 2 anti-automorphism: $\theta^2=\id_R$ and $\theta(ab)=\theta(b)\theta(a)$ for all $a,b\in R$.} $\theta$. For a right $R$-module $H$, define the left $R$-module $H^\theta$ by setting $H^\theta=H$ as abelian groups and $r\cdot h:=h\theta(r)$ for all $r\in R$ and $h\in H$. So the involution $\theta$ allows one to rewrite any right $R$-module $H$ as a left $R$-module $H^\theta$.

Recall that $\Hom_R(-,R)$ is a contravariant functor from the category of left $R$-modules to the category of right $R$-modules. We use the involution $\theta$ to transform the resulting right $R$-module into a left $R$-module. We call the resulting functor $\Hom_R^\theta(-,R)$. The idea is to reduce the computation of the homomorphism module to computing a kernel. Since $\Hom_R(-,R)$ is left exact, one obtains the exact sequence of right $R$-modules:
\[
  \xymatrix@1{
    0 \ar[r] & \Hom_R(M,R) \ar[r]^{\!\!\!\nu_M^*} & \Hom_R(R^{1\times l_0},R)
    \ar[r]^{\tM\cdot} & \Hom_R(R^{1\times l_1},R).
  }
\]
We compute $\Hom_R(M,R)$ as the kernel of the right most morphism. Following the row convention we identify $\Hom_R(R^{1\times j},R)$ with $R^{j\times 1}$, which justifies the notation $(\tM\cdot)$ for the right most morphism. Applying $\theta$ yields the exact sequence of left $R$-modules:
\[
  \xymatrix@1{
    0 \ar[r] & \Hom_R^\theta(M,R) \ar[r]^(.65){\nu_M^*} & R^{1\times l_0}
    \ar[r]^{\tM^\theta} & R^{1\times l_1},
  }
\]
where $(\tM^\theta)_{ab}=\theta(\tM_{ba})$ and $(R^{l\times 1})^\theta=R^{1\times l}$. Thus set
\[
  \ttHomR(M):=\Kernel(R^{1\times l_0} \xrightarrow{\tM^\theta} R^{1\times l_1}).
\]

This definition proposes setting $\Hull_\ttHomR(M):=R^{1\times l_0}$ and taking the embedding of the kernel $\ttHomR(M)$ in $R^{1\times l_0}$ as the natural embedding (cf.~\ref{ker}). The morphism part of the hull functor is hence defined by $\Hull_{\ttHomR}(\phi):=\phi^\theta$:
\[
  \xymatrix@R=1.5pc{
    0 \ar[r] & \ttHomR(M) \ar[r] & R^{1\times l_0} = (R^{l_0\times 1})^\theta \\
    0 \ar[r] & \ttHomR(N) \ar[u]^{\ttHomMapR(\phi)} \ar[r] &
    R^{1\times l_0''} = (R^{l_0''\times 1})^\theta \ar[u]_{\phi^\theta},
  }
\]
where $\phi^{\theta}:R^{1\times l_0}\to R^{1\times l_0''}$ is a morphism of free {\em left} modules.

Although the definition of $\Hull_\ttHomR(M)$ depends on the presentation of $M$, it is nevertheless functorial because of the equivalence of the categories $\CC^\mathrm{mod}$ and $\CC^\mathrm{mat}$.

$\ttHomR$ is a contravariant functor. The transposition in the definition of $\tM^\theta$ (and $\phi^\theta$) is the manifestation of this contravariance.

\subsubsection{The functor $\TensorProduct$}\label{tp}
Let $R$ be a {\em commutative}\footnote{Cf.\ Subsection~\ref{comm}.} ring. Recall that $-\otimes_R-$ is a bifunctor, covariant in both arguments. The idea is to reduce the computation of the tensor product module to computing a cokernel. Since $-\otimes_R-$ is right exact in both arguments, the tensor product of the two presentations $R^{1\times l_1}\xrightarrow{\tM} R^{1\times l_0}$ of $M$ and $R^{1\times l_1'}\xrightarrow{\tL} R^{1\times l_0'}$ of $L$ (each regarded as a  two term complex) is a presentation of $M\otimes_R L$ (cf.~\cite[Section~2.7]{GP}):
\[
  \xymatrix{
    (R^{1\times l_1}\otimes_R R^{1\times l_0'}) \oplus_R
    (R^{1\times l_0}\otimes_R R^{1\times l_1'}) \ar[r]^(.67)*+[o][F-]{\tau} &
    R^{1\times l_0}\otimes_R R^{1\times l_0'} \ar[rr]^(.55){\nu_{M\otimes_R L}} &&
    M\otimes_R L \ar[r] & 0.
  }
\]
We compute $M\otimes_R L$ as the cokernel of $\tau$. After identifying $R^{1\times j}\otimes_R R^{1\times k}$ with $R^{1\times jk}$ the morphism $\tau$ is given by the matrix $\mathtt{T}=\left(\begin{array}{c} \tM\otimes I_{l_0'} \\ I_{l_0}\otimes \tL \end{array}\right)$, where $\otimes$ is the {\sc Kronecker} product of matrices and $I_{l_0}$ resp.\  $I_{l_0'}$ is the identity matrix on $M$ resp.\  $L$.

Our convention in defining the {\sc Kronecker} product of two matrices $A=\left(a_{ij}\right)$ and $B$ is the usual one: $A\otimes B:=(a_{ij}B)$.

We define
\[
  \TensorProduct(M,L):=\Cokernel(R^{1\times (l_1l_0'+l_0l_1')}
  \xrightarrow{\mathtt{T}} R^{1\times l_0l_0'}).
\]

Define the morphism part of the functor $\TensorProduct(-,L)$ to be the {\sc Kronecker} product
\[
  \TensorProduct(\phi,L):=\phi\otimes I_{l_0'}.
\]
Here the hull functor $\Hull_{\TensorProduct(-,L)}$ coincides with the functor and the natural embedding is the identity transformation.

Analogously for the functor $\TensorProduct(L,-)$ define
\[
  \TensorProduct(L,\phi):=I_{l_0'}\otimes\phi.
\]
Again, the hull functor $\Hull_{\TensorProduct(-,L)}$ coincides with the functor and the natural embedding is the identity transformation.

$\TensorProduct(-,-)$ is a bifunctor, covariant in each argument.

\subsubsection{The functor $\ttHom$}\label{hom}
Let $R$ be a {\em commutative}\footnote{Cf.\ Subsection~\ref{comm}.} ring. Recall that $\Hom_R(-,-)$ is a bifunctor, contravariant in its first argument and covariant in the second. The idea is to reduce the computation of the homomorphism module to computing a kernel. Since $\Hom_R(-,L)$ is left exact for any module $L$ and $\Hom_R(P,-)$ is exact for $P$ projective (or free) one obtains (cf.~\cite[p.~104]{GP}):
\[
  \xymatrix@R=1.5pc{
    & & 0 & 0 \\
    0 \ar[r] & \Hom_R(M,L) \ar[r]^{\nu_M^*} & \Hom_R(R^{1\times l_0},L) \ar[u]
    \ar[r]^*+[o][F-]{\kappa} & \Hom_R(R^{1\times l_1},L) \ar[u] \\
    0 \ar[r] & \Hom_R(M,R^{1\times l_0'}) \ar[u]^{\nu_L} \ar[r]^{\nu_M^*} &
    \Hom_R(R^{1\times l_0},R^{1\times l_0'}) \ar[r]^{\tM\cdot} \ar[u]^{\nu_L} &
    \Hom_R(R^{1\times l_1},R^{1\times l_0'}) \ar[u]^{\nu_L} \\
    0 \ar[r] & \Hom_R(M,R^{1\times l_1'}) \ar[r]^{\nu_M^*} \ar[u] &
    \Hom_R(R^{1\times l_0},R^{1\times l_1'}) \ar[r]^{\tM\cdot} \ar[u]^{\cdot\tL} &
    \Hom_R(R^{1\times l_1},R^{1\times l_1'}) \ar[u]^{\cdot\tL}.
  }
\]
We compute $\Hom_R(M,L)$ as the kernel of $\kappa$ in the first row. Identifying $\Hom_R(R^{1\times j},R^{1\times k})$ with $R^{j\times k}$ justifies the notation used for the morphisms of the lower right square of the above diagram. Further, identifying $R^{j\times k}$ with $R^{1\times jk}$ (by writing all the $j$ rows as one long row) gives rise to the identification\footnote{$R^{j\times l_1'}\xrightarrow{\cdot\tL}R^{j\times l_0'}$ becomes $R^{1\times j l_1'}\xrightarrow{I_{j}\otimes \tL} R^{1\times jl_0'}$, and $R^{l_0\times k}\xrightarrow{\tM\cdot}R^{l_1\times k}$ becomes $R^{1\times l_0 k}\xrightarrow{(\tM\otimes I_{k})^\mathtt{tr}} R^{1\times l_1 k}$.} of $\Hom_R(R^{1\times l_0},L)$ with $\Cokernel(I_{l_0}\otimes \tL)$ and of $\Hom_R(R^{1\times l_1},L)$ with $\Cokernel(I_{l_1}\otimes \tL)$. The induced morphism
\[
  \kappa:\Hom_R(R^{1\times l_0},L)\to\Hom_R(R^{1\times l_1},L)
\]
is then given by the matrix $(\tM\otimes I_{l_0'})^\mathtt{tr}$. Thus define
\[
  \ttHom(M,L):=\Kernel(\Cokernel(I_{l_0}\otimes \tL) \xrightarrow{(\tM\otimes I_{l_0'})^\mathtt{tr}} \Cokernel(I_{l_1}\otimes \tL)).
\]

This definition proposes setting $\Hull_{\ttHom(-,L)}(M):=\Cokernel(I_{l_0}\otimes \tL)$ and taking the embedding of the kernel $\ttHom(M,L)$ in $\Cokernel(I_{l_0}\otimes \tL)$ as the natural embedding (cf.~\ref{ker}). The morphism part of the hull functor is hence defined by $\Hull_{\ttHom(-,L)}(\phi):=(\phi\otimes I_{l_0'})^\mathrm{tr}$:
\[
  \xymatrix@R=1.5pc{
    0 \ar[r] & \ttHom(M,L) \ar[r] &  \Cokernel(I_{l_0}\otimes \tL) \\
    0 \ar[r] & \ttHom(N,L) \ar[u]^{\ttHomMap(\phi)} \ar[r] &
    \Cokernel(I_{l_0''}\otimes \tL) \ar[u]_{(\phi\otimes I_{l_0'})^\mathrm{tr}}.
  }
\]

To address the functoriality of $\ttHom(-,-)$ in the second argument, interchange the role of $M$ and $L$, and set $\Hull_{\ttHom(L,-)}(M):=\Cokernel(I_{l_0'}\otimes \tM)$ and $\Hull_{\ttHom(L,-)}(\phi):=I_{l_0'}\otimes\phi$:
\[
  \xymatrix@R=1.5pc{
    0 \ar[r] & \ttHom(M,L) \ar[r] \ar[d]_{\mathtt{Hom2Map}(\phi)}  &
    \Cokernel(I_{l_0'}\otimes \tM) \ar[d]^{I_{l_0'}\otimes\phi} \\
    0 \ar[r] & \ttHom(N,L) \ar[r] & \Cokernel(I_{l_0'}\otimes \tN).
  }
\]

$\ttHom(-,-)$ is a bifunctor, contravariant in its first and covariant in its second argument. Transposing the matrix $\tM\otimes I_{l_0'}$ (and $\phi\otimes I_{l_0'}$) is the manifestation of the contravariance in the first argument.

\subsubsection{One last word on the commutativity of the ring $R$.} \label{comm}
For the most general definition of tensor product of modules one starts with a not necessarily commutative ring $R$, a right $R$-module $M_R$, and a left $R$-module ${}_R N$. If the $R$-module structure of $M$ resp.\  $N$ comes from a $(Q,R)$-bimodule resp.\  an $(R,S)$-bimodule structure, then their tensor product over $R$ is in a natural way a $(Q,S)$-bimodule ${}_Q M_R \otimes_R {}_R N_S$. $Q$ and $S$, again, are not necessarily commutative rings. Analogously, let $M$ and $N$ be two left modules over a not necessarily commutative ring $R$ and denote by $\Hom_R({}_R M, {}_R N)$ the abelian group of $R$-homomorphisms. If the $R$-module structure of $M$ resp.\  $N$ comes from an $(R,Q)$-bimodule resp.\  an $(R,S)$-bimodule structure, then $\Hom_R({}_R M, {}_R N)$ is again in a natural way a $(Q,S)$-bimodule.

Note that in both cases the resulting module might not be finitely generated as a $(Q,S)$-bimodule, even if $M$ and $N$ are finitely generated as $R$-modules.

In the special case where either $M$ or $N$ is an $(R,R)$-bimodule, then $M\otimes_R N$ resp.\ $\Hom_R(M,N)$ is again an $R$-module\footnote{In theory, this special case could work in $\homalg$, and indeed it does for very special cases, of course, beside the trivial case when $M=R$ or $N=R$.}. This is always the case when the ring is commutative, since then every $R$-module is an $(R,R)$-bimodule in the obvious way.

The {\sf Maple} implementation of $\homalg$ does neither support changing the ring, nor bimodule structures. These are issues we want to address in future implementations. But even though, the above mentioned problem of non-finite generation will remain the major obstacle.

\section{Derived functors} \label{Derivation}

The philosophy of derived categories is, roughly speaking, to replace a module by one of its resolutions, and then to look at the resolution as a very special type of complexes, with homology concentrated at degree $0$. After inverting quasi-isomorphisms, one obtains the derived category, where the objects are quasi-isomorphism types of complexes. Especially, all resolutions of a module become isomorphic objects in the derived category.

Using the cylinder-cone-translation construction \cite[{III.3}]{GM} one constructs out of every short exact sequence of complexes, a so called {\em distinguished triangle}. By passing to homology we again obtain distinguished triangles in the category of graded objects (which one can view as cyclic complexes, i.e.\ complexes with zero boundary maps \cite[III.2.3]{GM}). A popular way to start, is to take a distinguished triangle coming from a short exact sequence of complexes that are simultaneously resolving a short exact sequence of modules. Then one applies a functor, that turns such distinguished triangles again into distinguished triangles, and at last one takes the homology. The classical way of writing such a distinguished triangle of homologies is as a {\em long exact homology sequence}.

The reason for recalling the standard definitions in the following subsections is not only to indicate how they are computed using $\homalg$, but to finish the discussion of Section~\ref{mor}. This is done in Subsection \ref{rm}. 

\subsection{The procedure {\tt ResolveModule}: Resolve a module}\label{df}
By resolving a module $M$, which we view as a complex concentrated in degree $0$, we obtain a complex of free (resp.\  projective) modules and a quasi-isomorphism
\[
  \xymatrix@C=1.5pc@R=0.8pc{
    \cdots\ar[r] & P_2 \ar[r]\ar[d] & P_1  \ar[r]\ar[d] & P_0 \ar[r]\ar[d] & 0 \\
    \cdots\ar[r] & 0\ar[r] & 0\ar[r] & M \ar[r] & 0.
  }
\]
After inverting quasi-isomorphisms all resolutions become isomorphic.

For an additive functor $F$ and a module $M$ with a resolution
\[
  P:\xymatrix@1{ P_{q+1} \ar[r] & P_q \ar[r] &P_{q-1}\ar[r] &\cdots \ar[r] & P_1  \ar[r] & P_0 \ar[r] & 0}
\]
define the $q^\mathrm{th}$ derived functor $\LL_qF$ applied to $M$ by setting $\LL_qF(M):=H_q(F(P))$, the defect of the two consecutive morphisms $F(P_{q+1})\to F(P_q)$ and $F(P_q)\to F(P_{q-1})$ (cf.\ \ref{def} {\tt DefectOfHoms}). 

There is a cheaper method to compute the left derivation of a {\em right exact} covariant functor, which is based on \cite[the definition of $\widetilde{\LL}_q$ in  IV.(10.1), p.~156 and Prop.~IV.5.5, p.~133]{HS} and uses $\Kernel$ instead of $\DefectOfHoms$. By duality, there is a cheaper method to compute the right derivation of a {\em left exact} contravariant functor, which is based on \cite[Prop.~IV.5.8]{HS} and uses $\Cokernel$ instead of $\DefectOfHoms$. Both methods are implemented in $\homalg$.

\subsection{The procedure {\tt ResolutionOfSeq}: Resolve a morphism}\label{rm}
Given a morphism $M\xrightarrow{\phi} N$ one resolves $M$ and $N$ freely
\begin{equation}\tag{Lift}\label{Lift}
  \xymatrix@C=1.5pc@R=0.8pc{
    P_{q+1}\ar[r] & P_q \ar[r]\ar[d]^{\phi_q} & P_{q-1}\ar[r]\ar[d]^{\phi_{q-1}} &
    \cdots \ar[r] & P_1  \ar[r]\ar[d]^{\phi_1} & P_0
    \ar[r]\ar[d]^{\phi_0} & M\ar[d]^\phi \ar[r] & 0 \\
    P_{q+1}'\ar[r] & P_q'\ar[r] &  P_{q-1}'\ar[r] & \cdots \ar[r] & P_1'\ar[r] & P_0' \ar[r] & N\ar[r] & 0,
  }
\end{equation}
and computes the $\phi_q$'s by iteratively completing image squares. Here one needs a free resolution of $M$ to be in case (2) of {\tt CompleteImSq}, \ref{CompleteImSq}. Applying $F$ one gets $F(P_q)\xrightarrow{F(\phi_q)} F(P_q')$. But now the object part of the functor $\LL_qF$ applied to $M$ (resp.\  $N$) is a subfactor module of $F(P_q)$ (resp.\  $F(P_q')$) and one is in the situation of $\DefectOfHoms$, \ref{def}. This is all what $\FunctorMap$ needs to compute the morphism part of a derived functor. This finishes the discussion of Section~\ref{mor}. 

\subsection{The procedure {\tt ResolveShortExactSeq}: Resolving a short exact sequence of modules}
To resolve a short exact sequence of modules $0\to M' \to M \to M'' \to 0$, one starts with a resolution of $M''$:
\[
  \xymatrix@C=1.5pc@R=0.8pc{
    & & & & 0\ar[d] & \\
    & & & & M'\ar[d]\ar[r] & 0\\
    & & & & M\ar[d]\ar[r] & 0 \\
    \cdots\ar[r] & P_2'' \ar[r]\ar[d] & P_1''  \ar[r]\ar[d] & P_0'' \ar[r]\ar[d] & M''  \ar[r]\ar[d] & 0. \\
    & 0 & 0 & 0 & 0 &
  }
\]
Then one completes the middle line by taking free hulls of iterated pullbacks. Finally one fills the upper line by taking kernels to obtain:
\begin{equation}\label{M}\tag{$M$}
  \xymatrix@C=1.5pc@R=0.8pc{
    & 0\ar[d] & 0\ar[d] & 0\ar[d] & 0\ar[d] & \\
    \cdots\ar[r] & P_2'\ar[r]\ar[d] & P_1'\ar[r]\ar[d] & P_0'\ar[r]\ar[d] & M'\ar[d]\ar[r] & 0\\
    \cdots\ar[r] & P_2\ar[r]\ar[d] & P_1\ar[r]\ar[d] & P_0\ar[r]\ar[d] & M\ar[d]\ar[r] & 0 \\
    \cdots\ar[r] & P_2'' \ar[r]\ar[d] & P_1''  \ar[r]\ar[d] & P_0'' \ar[r]\ar[d] & M''  \ar[r]\ar[d] & 0, \\
    & 0 & 0 & 0 & 0 &
  }
\end{equation}
with exact columns and rows. This method is implemented in the procedure {\tt ResolveShort\-ExactSeq}. There are several other methods to resolve a short exact sequence simultaneously, cf.\  \cite[Proof of Theorem IV.6.1]{HS}.

Computing the pullback $A'$ of $B'\xrightarrow{\beta}B\xleftarrow{\phi}A$
is reduced to computing a kernel, namely that of $A\oplus B'\xrightarrow{\left(\begin{array}{c}\phi \\ -\beta\end{array}\right)}B$. The two maps $B'\xleftarrow{\psi}A'\xrightarrow{\alpha}A$ of the pullback are the two parts of the kernel embedding $0\xrightarrow{} A'\xrightarrow{\left(\begin{array}{cc}\alpha&\psi\end{array}\right)}A\oplus B'$ (cf.~\ref{ker}).

\subsection{The procedure {\tt LongExactHomologySeq}: Connecting homomorphism and long exact sequences}\label{ch}
Applying an additive covariant functor $F$ to the truncation of the diagram (\ref{M})
\begin{equation}\label{P}\tag{$P$}
  \xymatrix@C=1.5pc@R=0.8pc{
    & 0\ar[d] & 0\ar[d] & 0\ar[d] & \\
    \cdots\ar[r] & P_2'\ar[r]\ar[d] & P_1'\ar[r]\ar[d] & P_0'\ar[r]\ar[d] & 0\\
    \cdots\ar[r] & P_2\ar[r]\ar[d] & P_1\ar[r]\ar[d] & P_0\ar[r]\ar[d] & 0 \\
    \cdots\ar[r] & P_2'' \ar[r]\ar[d] & P_1''  \ar[r]\ar[d] & P_0'' \ar[r]\ar[d] & 0 \\
    & 0 & 0 & 0 &
  }
\end{equation}
results in a diagram $\overline{P}=F(P)$ where the columns are still exact, but where the rows are now in positive degrees, in general, no longer exact. Connecting the homologies of the rows one obtains the {\em long exact homology sequence}

\[
  \cdots \xrightarrow{\partial_{q+1}} \LL_qF(M')\to \LL_qF(M)\to \LL_qF(M'')\xrightarrow{\partial_q} \LL_{q-1}F(M')\to\cdots,
\]
where the {\em connecting homomorphisms} are computed via the snake lemma applied to the diagrams
\[
  \xymatrix@C=1.5pc@R=0.8pc{
    & 0\ar[d] \\
    \overline{P}_{q+1}'/\overline{B}_{q+1}'\ar[r]\ar[d] & \overline{Z}_q' \ar[d] \\
    \overline{P}_{q+1}/\overline{B}_{q+1}\ar[r]\ar[d] & \overline{Z}_q\ar[d] \\
    \overline{P}_{q+1}''/\overline{B}_{q+1}'' \ar[r]\ar[d] & \overline{Z}_q'', \\
    0 &
  }
\]
by taking kernels and cokernels, where $\overline{B}_i:=\img(\overline{P}_{i+1}\to \overline{P}_i)$ and $\overline{Z}_i:=\ker(\overline{P}_i\to \overline{P}_{i-1})$. One computes the connecting homomorphism from $\ker(\overline{P}_{q+1}''/\overline{B}_{q+1}'' \to \overline{Z}_q'')$ to $\coker( \overline{P}_{q+1}'/\overline{B}_{q+1}'\to \overline{Z}_q' )$ by a diagram chase, which accounts in taking preimages twice, namely under the maps $\overline{P}_{q+1}/\overline{B}_{q+1}\to \overline{P}_{q+1}''/\overline{B}_{q+1}''$ and $ \overline{Z}_q'\to \overline{Z}_q$. For this the procedure {\tt Preimage} from \ref{pi} is used.

Applying a contravariant functor to (\ref{P}) yields a {\em long exact cohomology sequence}. The corresponding procedure is called {\tt LongExactCohomologySeq}.

\appendix

\section{The ring packages}\label{RP}
The following {\sf Maple} ring packages have successfully been used with $\homalg$. In each of the following descriptions we append a list of rings which can be dealt with in $\homalg$ using the respective package. Further, and without any extra help from the ring package, $\homalg$ can automatically compute over residue class rings of any supported ring.
\begin{itemize}
\item {\tt PIR} \cite{PIR} is one more tiny package, or rather a pseudo-package, that makes {\sf Maple}'s built-in facilities for dealing with integers and some other principal ideal rings available to $\homalg$.
(Prime subfields $\Q$ and $\Z/p\Z$ and their finite field extensions, realized as primitive extensions, rational function fields over the previous fields, the integers $\Z$, the {\sc Gauss}ian integers $\Z[\sqrt{-1}]$ and univariate polynomial rings $\Z/p\Z[x]$, where $p$ is a prime, $\Q[x]$ and $K[x]$, where $K$ is a rational function field over a finite extension of $\Q$, realized as a primitive extension.)
\item {\tt Involutive}~\cite{JanetPackage} implements the involutive basis
  technique of {\sc V.~P.\ Gerdt} and {\sc Y.~A.\ Blinkov} in {\sf Maple}. An
  involutive basis is a special kind of {\sc Gr{\"o}bner} basis for an ideal
  of a polynomial ring or, more generally, for a submodule of a free module
  over a polynomial ring. Involutive bases have nice combinatorial properties
  \cite{PRJ}, and the algorithms designed by {\sc V.~P.~Gerdt} and {\sc
    Y.~A.~Blinkov} \cite{GerI,GB1,GB2}  compute them efficiently. In fact,
  these algorithms provide an efficient alternative to {\sc Buchberger}'s
  algorithm \cite{Buch} to compute {\sc Gr{\"o}bner} bases. {\tt Involutive}
  restricts to particular involutive bases, namely {\sc Janet} bases. It also
  provides an interface to a {\sf C++} implementation of the involutive basis
  technique which can be used to call the fast routines when needed as well as
  to switch to these fast routines for the whole {\sf Maple} session.
  (Commutative polynomial rings: $S[x_1, \ldots, x_n]$, where $S$ is either $\Z$ or a field existing in {\sf Maple}.)
\item {\tt Janet}~\cite{JanetPackage} implements the involutive basis technique for computing
 {\sc Janet} bases of linear systems of partial differential equations. (Differential algebras over differential fields:
$K[\frac{\partial}{\partial x_1}, \ldots, \frac{\partial}{\partial x_n}]$, where $K$ is a differential field which exists in {\sf Maple}.)
\item {\tt JanetOre}~\cite{robphd,JO} generalizes {\tt Involutive} from commutative polynomial rings to certain iterated skew polynomial rings. In particular, it computes Janet bases for left ideals in Ore algebras \cite{CSO}. ($K[\partial; \sigma, \delta]$, where $K$ is a polynomial ring over a field, $\partial$ a new indeterminate, $\sigma$ is a certain automorphism of $K$ and $\delta$ a $\sigma$-derivation of $K$, and iterated extensions of this kind.)
\item {\tt OreModules}~\cite{CQR07} is a Maple package for the study of structural properties of linear systems over {\sc Ore} algebras, i.e.\ linear equations involving certain linear functional operators which can be considered as elements of an {\sc Ore} algebra. By default, it uses the {\sf Maple} package {\tt Ore\_{}algebra} \cite{CSO} to compute {\sc Gr{\"o}bner} bases, but these calls can also be switched to {\tt JanetOre}. ({\sc Ore} algebras \cite{CSO} and the iterated skew polynomial rings from the previous point.)
\end{itemize}
$\homalg$ is also able to make use of various normal form algorithms for modules resp.\   special types of modules over various rings, which are used to provide a standard form for a presentation of these modules:
\begin{itemize}
  \item {\tt PIR} uses the {\sc Smith} normal form for
    ($\mathsf{Maple}$-built-in) principal ideal rings.
  \item {\tt Janet} optionally uses the {\sc Jacobson} normal form for
    univariate differential rings, i.e.\  rings of the form $K[\partial]$,
    where $K$ is a differential field with $\partial$ a derivation of $K$.
  \item {\tt Involutive} optionally uses the extension package
    {\tt QuillenSuslin} written by {\sc Anna Fabianska} \cite{Fab,FQ}
    implementing the {\sc Quillen-Suslin} theorem to
    compute a free basis of a projective module over a polynomial ring (which
    is then free by the theorem). A similar extension package is planned for
    {\tt OreModules}.
  \item {\tt OreModules} optionally uses the extension package {\tt Stafford} \cite{QR_S}
    which computes a free basis for a stably free module of rank at least $2$ over the Weyl
    algebras $k[x_1,\ldots,x_n,\partial_1,\ldots,\partial_n]$ and
    $k(x_1,\ldots,x_n)[\partial_1,\ldots,\partial_n]$, with $k$ a field of characteristic
    $0$.
\end{itemize}

\section{The {\sf Maple} implementation}\label{implementation}

\subsection{Presentations}
A presentation of a module in the current $\mathsf{Maple}$ implementation of $\homalg$ is a list\footnote{In future implementations we will make use of object oriented data structures, which encapsulate all this information.} containing as first entry the list of generators and as second entry the list of relations. The third entry is a string delimiter to optically indicate the end of the presentation. This string, unless changed by the user, defaults to {\tt "Presentation"}. The remaining entries provide extra information about the presented module, e.g.\ its {\sc Hilbert} series in case the ring is the polynomial ring. This extra information can only be provided by the ring-specific package.

For $M=\coker(\tM)=R^{1\times l_0}/R^{1\times l_1}\tM$ the list of the concrete generators are numbered by abstract generators being the $l_0$ standard basis row vectors of the underlying free module $R^{1\times l_0}$. The list of relations simply contains the rows of the matrix $\tM$. An illustration is given in Figure \ref{pres} of Example~\ref{ExHom}.

\subsection{The morphism part of functors}\label{Basic}
The name convention for functors used in the {\sf Maple} implementation of $\homalg$ is as follows: If the procedure implementing the object part of a functor is called {\tt F}, then the procedure implementing the morphism part is called {\tt FMap}. It is defined using the procedure $\FunctorMap$ applied to the object part procedure {\tt F}.
The several pieces of code for the morphism part of the functors \href{http://wwwb.math.rwth-aachen.de/homalg/implementation.html#CokernelMap}{$\Cokernel$},
\href{http://wwwb.math.rwth-aachen.de/homalg/implementation.html#KernelMap}{$\Kernel$},
\href{http://wwwb.math.rwth-aachen.de/homalg/implementation.html#DefectOfHomsMap}{$\DefectOfHoms$} and \href{http://wwwb.math.rwth-aachen.de/homalg/implementation.html#HomMap_R}{$\ttHomR$},
and the two bifunctors \href{http://wwwb.math.rwth-aachen.de/homalg/implementation.html#TensorProductMap}{$\TensorProduct$} and \href{http://wwwb.math.rwth-aachen.de/homalg/implementation.html#Hom}{$\ttHom$},
reproduced on the web \cite{homalg_implementation}, demonstrate how $\FunctorMap$ unifies the definition of the morphism part of {\em all} functors in $\homalg$.

In future implementations of $\homalg$ the two procedures implementing the object and morphism part of a functor will be unified in one. The unified procedure will be able to recognize if it has been applied to an object, to a morphism or even to complexes. This will be an easy task, once we strictly use structures throughout $\homalg$.

\subsection{Encapsulating functors}
A functor is fully defined when both parts are defined, i.e.\  its object part and its morphism part. If the functor is a multi-functor, then several morphism parts have to be defined. $\homalg$ accesses all these parts of a functor via a so called encapsulation. It incorporates both parts and also takes care of the possible multi-functoriality. The implementation of the encapsulations of the functors \href{http://wwwb.math.rwth-aachen.de/homalg/implementation.html#Functor_Hom_R}{$\ttHomR$}, \href{http://wwwb.math.rwth-aachen.de/homalg/implementation.html#Functor_TensorProduct}{$\TensorProduct$} and \href{http://wwwb.math.rwth-aachen.de/homalg/implementation.html#Functor_Hom}{$\ttHom$} can be viewed under \cite{homalg_implementation}.

\subsection{Composition of functors}\label{ComposeFunctors}
As was mentioned in Section~\ref{mor}, composing two functors is an easy task since one simply has to compose their actions on objects and on morphisms. The procedure responsible for this is called {\tt ComposeFunctors}. The implementation of the composition of the functors 
\href{http://wwwb.math.rwth-aachen.de/homalg/implementation.html#HomHom_R}{$\ttHomHomR$} and \href{http://wwwb.math.rwth-aachen.de/homalg/implementation.html#HomHom}{$\ttHomHom$} can be viewed under \cite{homalg_implementation}.

\subsection{Applying functors to complexes}
As mentioned in Section~\ref{Derivation}, functors should be applied to complexes of modules, rather than to single modules. Thereby $\homalg$ again makes use of the encapsulation of functors. In this implementation, there are two procedures depending on whether one is dealing with a covariant or a contravariant functor. They are called {\tt FunctorOnSeqs} and {\tt CofunctorOnSeqs}. The implementation of the functors \href{http://wwwb.math.rwth-aachen.de/homalg/implementation.html#HomOnSeqs}{$\ttHom$} and \href{http://wwwb.math.rwth-aachen.de/homalg/implementation.html#HomHomOnSeqs}{$\ttHomHom$} on complexes can be viewed under \cite{homalg_implementation}.

\subsection{Derivation of functors}\label{der}
The left derivation procedure for covariant functors is called {\tt LeftDerivedFunctor}. The faster derivation procedure for right exact functors, for example $-\otimes L$, referred to at the end of Subsection~\ref{df} is called {\tt LeftDerived\-RightExactFunctor}.

The right derivation procedure for contravariant functors is called {\tt RightDerived\-Co\-functor}. The faster derivation procedure for left exact contravariant functors, for example $\Hom(-,L)$, referred to at the end of Subsection~\ref{df} is called {\tt RightDerived\-LeftExact\-Cofunctor}.

The implementation of the left derived functor \href{http://wwwb.math.rwth-aachen.de/homalg/implementation.html#LHomHom}{$\ttLHomHom$} (of $\ttHomHom$ with respect to its first argument) and the right derived functor \href{http://wwwb.math.rwth-aachen.de/homalg/implementation.html#Ext}{$\ttExt$} (of $\ttHom$ with respect to its first argument) can be viewed under \cite{homalg_implementation}.

\subsection{How $\FunctorMap$ works}\label{FunctorMap}
Let {\tt F} denote the object part procedure of a functor $F$. First $\FunctorMap$ asks {\tt F} if the underlying functor is co- or contravariant. If {\tt F} is $\homalg$-basic it has to know the answer itself. If {\tt F} is defined as a composition $F_1\circ F_2$ then the question is passed to the procedure {\tt ComposeFunctors}, which decides the answer by asking $F_1$ and $F_2$ (this is recursive). If {\tt F} is defined as a derivation $\LL_q G$ or $\RR^q G$, using one of the derivation procedures, then {\tt F} passes the question to the latter, which decides the answer by asking $G$ (this is recursive). Any recursion ends when a $\homalg$-basic functor is reached.

$\FunctorMap$ then asks {\tt F} if it is defined using the procedure {\tt ComposeFunctors}. If {\tt F} is $\homalg$-basic it ignores the question. If {\tt F} is defined using one of the derivation procedures, the question is passed to the latter, which ignores it. If {\tt F} is indeed defined as a composition $F_1\circ F_2$, then {\tt F} passes the question to {\tt ComposeFunctors}, which returns the two functors $F_1$ and $F_2$ in both their parts. $\FunctorMap$ can now easily construct the morphism part of $F$ by composing the morphism parts of $F_1$ and $F_2$. In case {\tt F} is not defined as a composition, $\FunctorMap$ asks it if is defined by derivation. If {\tt F} is $\homalg$-basic it ignores the question. If {\tt F} is indeed defined as a derivation $\LL_q G$ or $\RR^q G$, using one of the derivation procedures, then {\tt F} passes the question to the latter, which returns the functor $G$ in both its parts {\em and} a procedure based on {\tt ResolutionOfSeq} to compute $\phi_q$ out of $\phi$ (cf.~(\ref{Lift}), p.~\pageref{Lift}). With these two ingredients $\FunctorMap$ is able to construct the morphism part of {\tt F} as described in Subsection \ref{rm}.  If {\tt F} is defined neither by composition nor by derivation, i.e.\  its is $\homalg$-basic, then it is asked by $\FunctorMap$ to return\footnote{The technical details to realize this heavily depend on the implementation.} its hull functor $\Hull_\mathtt{F}$ together with the natural embedding (cf.~(\ref{Hull}), p.~\pageref{Hull}). Again this suffices to construct the morphism part of {\tt F} as described in Section \ref{mor}.

\section{Examples}\label{Ex}

For further examples we refer to \cite{BRMTNS,BREACA} and the site of $\homalg$ \cite{homalg}.

A nice application is the tiny $\homalg$-based package $\mathtt{conley}$ that computes $C$-connection matrices of graded module octahedra/braids of Morse decompositions in dynamical system theory \cite{BRconley}.

For the sake of demonstration we wrote a tiny package called {\tt alexander} \cite{alexander}, which relies on $\homalg$ and computes simplicial homology and cohomology. Future implementations of $\homalg$ are planned to enable more serious applications to topology.

\subsection{Example 1}\label{ExHom}

\PUSH{Hom.tex}%
\DefineParaStyle{Maple Output}
\DefineCharStyle{2D Math}
\DefineCharStyle{2D Output}
\DefineCharStyle{LaTeX}
\begin{maplegroup}
In this example we compute a module of homomorphisms.

\end{maplegroup}
\begin{maplegroup}
\begin{mapleinput}
\mapleinline{active}{1d}{restart:}{%
}
\end{mapleinput}

\end{maplegroup}
\begin{maplegroup}
\begin{mapleinput}
\mapleinline{active}{1d}{with(Involutive): with(homalg):}{%
}
\end{mapleinput}

\end{maplegroup}
\begin{maplegroup}
Specify the {\tt homalg}-table of the ring package {\tt Involutive}:

\end{maplegroup}
\begin{maplegroup}
\begin{mapleinput}
\mapleinline{active}{1d}{RPI:=`Involutive/homalg`;}{%
}
\end{mapleinput}

\mapleresult
\begin{maplelatex}
\mapleinline{inert}{2d}{RPI := `Involutive/homalg`;}{%
\[
\mathit{RPI} := \mathit{Involutive/homalg}
\]
}
\end{maplelatex}

\end{maplegroup}
\begin{maplegroup}
Use the ring package {\tt Involutive} as the default ring package:

\end{maplegroup}
\begin{maplegroup}
\begin{mapleinput}
\mapleinline{active}{1d}{`homalg/default`:=RPI;}{%
}
\end{mapleinput}

\mapleresult
\begin{maplelatex}
\mapleinline{inert}{2d}{`homalg/default` := `Involutive/homalg`;}{%
\[
\mathit{homalg/default} := \mathit{Involutive/homalg}
\]
}
\end{maplelatex}

\end{maplegroup}
\begin{maplegroup}
Define the ring $R=\Q[x,y,z]$:

\end{maplegroup}
\begin{maplegroup}
\begin{mapleinput}
\mapleinline{active}{1d}{var:=[x,y,z];}{%
}
\end{mapleinput}

\mapleresult
\begin{maplelatex}
\mapleinline{inert}{2d}{var := [x, y, z];}{%
\[
\mathit{var} := [x, \,y, \,z]
\]
}
\end{maplelatex}

\end{maplegroup}
\begin{maplegroup}
\begin{mapleinput}
\mapleinline{active}{1d}{K:=Cokernel([[x,y,0],[x^2,y^2,0],[x^3,y^3,z^3]],var);}{%
}
\end{mapleinput}

\mapleresult
\begin{maplelatex}
\mapleinline{inert}{2d}{K := [[[1, 0, 0] = [1, 0, 0], [0, 1, 0] = [0, 1, 0], [0, 0, 1] = [0,
0, 1]], [[x, y, 0], [0, x*y-y^2, 0], [0, 0, z^3]], "Presentation",
3+8*s+14*s^2+s^3*(14/(1-s)+6/(1-s)^2), [14, 6, 0]];}{%
\maplemultiline{
K := [[[1, \,0, \,0]=[1, \,0, \,0], \,[0, \,1, \,0]=[0, \,1, \,0]
, \,[0, \,0, \,1]=[0, \,0, \,1]],  \\
[[x, \,y, \,0], \,[0, \,x\,y - y^{2}, \,0], \,[0, \,0, \,z^{3}]]
, \,\mbox{``Presentation''},  \\
3 + 8\,s + 14\,s^{2} + s^{3}\,({\displaystyle \frac {14}{1 - s}} 
 + {\displaystyle \frac {6}{(1 - s)^{2}}} ), \,[14, \,6, \,0]] }
}
\end{maplelatex}

\end{maplegroup}
\begin{maplegroup}
\begin{mapleinput}
\mapleinline{active}{1d}{L:=Cokernel([[x,y]],var);}{%
}
\end{mapleinput}

\mapleresult
\begin{maplelatex}
\mapleinline{inert}{2d}{L := [[[1, 0] = [1, 0], [0, 1] = [0, 1]], [[x, y]], "Presentation",
2+s*(2/(1-s)+2/(1-s)^2+1/((1-s)^3)), [2, 2, 1]];}{%
\maplemultiline{
L := [[[1, \,0]=[1, \,0], \,[0, \,1]=[0, \,1]], \,[[x, \,y]], \,
\mbox{``Presentation''},  \\
2 + s\,({\displaystyle \frac {2}{1 - s}}  + {\displaystyle 
\frac {2}{(1 - s)^{2}}}  + {\displaystyle \frac {1}{(1 - s)^{3}}
} ), \,[2, \,2, \,1]] }
}
\end{maplelatex}

\end{maplegroup}
\begin{maplegroup}
Compute the module of homomorphisms $\Hom_R(L,K)$ (see
Figure~\ref{pres}):

\end{maplegroup}
\begin{maplegroup}
\begin{mapleinput}
\mapleinline{active}{1d}{hom:=Hom(L,K,var);}{%
}
\end{mapleinput}

\mapleresult
\begin{maplelatex}
\mapleinline{inert}{2d}{hom := [[[1, 0, 0] = matrix([[1, 0, 0], [0, 1, 0]]), [0, 1, 0] =
matrix([[0, 0, -y], [0, 0, x]]), [0, 0, 1] = matrix([[0, 0, 0], [0,
-y+x, 0]])], [[0, 0, y], [x^2-x*y, 0, -x], [0, z^3, 0]],
"Presentation", 3+8*s+14*s^2+s^3*(14/(1-s)+6/(1-s)^2), [14, 6, 0]];}{%
\maplemultiline{
\mathit{hom} := [[[1, \,0, \,0]= \left[ 
{\begin{array}{rrr}
1 & 0 & 0 \\
0 & 1 & 0
\end{array}}
 \right] , \,[0, \,1, \,0]= \left[ 
{\begin{array}{rrc}
0 & 0 &  - y \\
0 & 0 & x
\end{array}}
 \right] , \,[0, \,0, \,1]= \left[ 
{\begin{array}{rcr}
0 & 0 & 0 \\
0 &  - y + x & 0
\end{array}}
 \right] ],  \\
[[0, \,0, \,y], \,[x^{2} - x\,y, \,0, \, - x], \,[0, \,z^{3}, \,0
]], \,\mbox{``Presentation''},  \\
3 + 8\,s + 14\,s^{2} + s^{3}\,({\displaystyle \frac {14}{1 - s}} 
 + {\displaystyle \frac {6}{(1 - s)^{2}}} ), \,[14, \,6, \,0]] }
}
\end{maplelatex}

\end{maplegroup}
\begin{maplegroup}
\begin{figure}[t]
\caption{\label{pres} A module of homomorphisms between two modules
over
$R=\Q[x,y,z]$ with {\tt Involutive}}

\begin{center}
\begin{picture}(380,20)
\put(0.0,0.0){
{\scriptsize
$[$
$[[1,0,0]= \left[ \begin {array}{rrr} 1 & 0 & 0
\\ \noalign{\medskip} 0 & 1 & 0
\end{array} \right],
[0,1,0]= \left[ \begin{array}{rrr} 0 & 0 & -y
\\ \noalign{\medskip} 0 & 0 & x \end{array} \right],
[0,0,1]= \left[ 
\begin{array}{rcr} 0 & 0 & 0 \\\noalign{\medskip} 0 & x-y & 0
\end{array}
 \right] ],$
}
}
\end{picture}
\end{center}

\bigskip

\begin{center}
\begin{picture}(260,10)
\put(0.0,0.0){\scriptsize
$[[0,0,y],[x^2-xy,0,-x],[0,z^{3},0]],$
}
\end{picture}
\end{center}

\begin{center}
\begin{picture}(180,10)
\put(0.0,0.0){
{\scriptsize ``Presentation''},
}
\end{picture}
\end{center}

\begin{center}
\begin{picture}(100,10)
\put(-100.0,10.0){\line(0,1){40}}
\put(-120.0,0.0){\scriptsize generators}
\put(0.0,0.0){\scriptsize
$3+8\,s+14\,s^{2}+s^{3}
\left( \frac{14}{(1-s)}+ \frac{6}{(1-s)^2} \right)$,
}
\end{picture}
\end{center}

\begin{center}
\hspace{1 cm}
\begin{picture}(60,10)
\put(-100.0,12.0){\line(0,1){30}}
\put(-120.0,0.0){\scriptsize relations}
\put(0.0,0.0){\scriptsize
$[14,6,0]$ $]$
}
\end{picture}
\end{center}

\begin{center}
\begin{picture}(60,10)
\put(-8.0,12.0){\line(0,1){15}}
\put(-45.0,0.0){\scriptsize {\sc Hilbert} series}
\end{picture}
\end{center}

\begin{center}
\begin{picture}(60,10)
\put(30.0,15.0){\line(0,1){12}}
\put(10.0,3.0){\scriptsize {\sc Cartan} characters}
\end{picture}
\end{center}

\end{figure}

\end{maplegroup}
%
\POP

\subsection{Example 2}

\PUSH{TwoNonIsomorphicExtensions.tex}%
\DefineParaStyle{Maple Output}
\DefineCharStyle{2D Math}
\DefineCharStyle{2D Output}
\DefineCharStyle{LaTeX}
\begin{maplegroup}
Let $R := \Q[x,y,z]$.
In this example we want to study two non-equivalent extensions $0 \to
K \to M \to L \to 0$ and $0 \to K \to N \to L \to 0$ with $K$ a
torsion module and $L$ a torsion free module. Our goal is to use the
notion of functor to reveal that these extensions are not only
non-equivalent but also non-isomorphic. For this we define two
functors $F_M$ and $F_N$ and study their behavior when applied to
complexes. Concretely, we apply $F_M$ resp.\  $F_N$ to $0 \to K \to M
\to L \to 0$ which we refer to as
\[ S: 0 \to K \xrightarrow{\alpha_1} M \xrightarrow{\alpha_2} L \to 0.
\]
In the following we consider
\[ F_Q=\ttExt(1,-,Q) = \Ext^1(-,Q) = R^1\Hom_R(-, Q) \]

\medskip

\end{maplegroup}
\begin{maplegroup}
\begin{mapleinput}
\mapleinline{active}{1d}{restart;}{%
}
\end{mapleinput}

\end{maplegroup}
\begin{maplegroup}
\begin{mapleinput}
\mapleinline{active}{1d}{with(Involutive): with(homalg):}{%
}
\end{mapleinput}

\end{maplegroup}
\begin{maplegroup}
Use the ring package $\mathtt{Involutive}$ as the default ring
package:

\end{maplegroup}
\begin{maplegroup}
\begin{mapleinput}
\mapleinline{active}{1d}{RPI:=`Involutive/homalg`;}{%
}
\end{mapleinput}

\mapleresult
\begin{maplelatex}
\mapleinline{inert}{2d}{RPI := `Involutive/homalg`;}{%
\[
\mathit{RPI} := \mathit{Involutive/homalg}
\]
}
\end{maplelatex}

\end{maplegroup}
\begin{maplegroup}
\begin{mapleinput}
\mapleinline{active}{1d}{`homalg/default`:=RPI;}{%
}
\end{mapleinput}

\mapleresult
\begin{maplelatex}
\mapleinline{inert}{2d}{`homalg/default` := `Involutive/homalg`;}{%
\[
\mathit{homalg/default} := \mathit{Involutive/homalg}
\]
}
\end{maplelatex}

\end{maplegroup}
\begin{maplegroup}
We force the $\mathsf{Maple}$ ring package $\mathtt{Involutive}$ to
use the external program $\mathtt{JB}$, which is a {\sf C++}
implementation of the involutive division algorithm:

\end{maplegroup}
\begin{maplegroup}
\begin{mapleinput}
\mapleinline{active}{1d}{InvolutiveOptions("C++");}{%
}
\end{mapleinput}

\end{maplegroup}
\begin{maplegroup}
Define the ring $R=\Q[x,y,z]$:

\end{maplegroup}
\begin{maplegroup}
\begin{mapleinput}
\mapleinline{active}{1d}{var:=[x,y,z];}{%
}
\end{mapleinput}

\mapleresult
\begin{maplelatex}
\mapleinline{inert}{2d}{var := [x, y, z];}{%
\[
\mathit{var} := [x, \,y, \,z]
\]
}
\end{maplelatex}

\end{maplegroup}
\begin{maplegroup}
The two presentation matrices $\mathtt{M}$ and $\mathtt{N}$:

\end{maplegroup}
\begin{maplegroup}
\begin{mapleinput}
\mapleinline{active}{1d}{MM:=matrix([[x*z, z*y, z^2, 0, 0, y], [0, 0, 0, z^2*y-z^2, z^3, x*z],
[0, 0, 0, z*y^2-z*y, z^2*y, x*y], [0, 0, 0, x*z*y-x*z, x*z^2, x^2],
[-x*y, -y^2, -z*y, x^2*y-x^2-y+1, x^2*z-z, 0], [x^2*y-x^2, x*y^2-x*y,
x*z*y-x*z, -y^3+2*y^2-y, -z*y^2+z*y, 0]]);}{%
}
\end{mapleinput}

\mapleresult
\begin{maplelatex}
\mapleinline{inert}{2d}{MM := matrix([[x*z, z*y, z^2, 0, 0, y], [0, 0, 0, z^2*y-z^2, z^3,
x*z], [0, 0, 0, z*y^2-z*y, z^2*y, x*y], [0, 0, 0, x*z*y-x*z, x*z^2,
x^2], [-x*y, -y^2, -z*y, x^2*y-x^2-y+1, x^2*z-z, 0], [x^2*y-x^2,
x*y^2-x*y, x*z*y-x*z, -y^3+2*y^2-y, -z*y^2+z*y, 0]]);}{%
\[
\mathit{MM} :=  \left[ 
{\begin{array}{cccccc}
x\,z & z\,y & z^{2} & 0 & 0 & y \\
0 & 0 & 0 & z^{2}\,y - z^{2} & z^{3} & x\,z \\
0 & 0 & 0 & z\,y^{2} - z\,y & z^{2}\,y & x\,y \\
0 & 0 & 0 & x\,z\,y - x\,z & x\,z^{2} & x^{2} \\
 - x\,y &  - y^{2} &  - z\,y & x^{2}\,y - x^{2} - y + 1 & x^{2}\,
z - z & 0 \\
x^{2}\,y - x^{2} & x\,y^{2} - x\,y & x\,z\,y - x\,z &  - y^{3} + 
2\,y^{2} - y &  - z\,y^{2} + z\,y & 0
\end{array}}
 \right] 
\]
}
\end{maplelatex}

\end{maplegroup}
\begin{maplegroup}
\begin{mapleinput}
\mapleinline{active}{1d}{NN:=matrix([[x*y, y^2, z*y, 0, 0, z], [0, 0, 0, z*y^2-z*y, z^2*y,
x*z], [x^2*z, x*z*y, x*z^2, -z^2*y+z^2, -z^3, 0], [0, 0, 0,
y^3-2*y^2+y, z*y^2-z*y, x*y-x], [0, 0, 0, x^2*y-x^2-y+1, x^2*z-z, y],
[x^3, x^2*y, x^2*z, -x*z*y+x*z, -x*z^2, 0]]);}{%
}
\end{mapleinput}

\mapleresult
\begin{maplelatex}
\mapleinline{inert}{2d}{NN := matrix([[x*y, y^2, z*y, 0, 0, z], [0, 0, 0, z*y^2-z*y, z^2*y,
x*z], [x^2*z, x*z*y, x*z^2, -z^2*y+z^2, -z^3, 0], [0, 0, 0,
y^3-2*y^2+y, z*y^2-z*y, x*y-x], [0, 0, 0, x^2*y-x^2-y+1, x^2*z-z, y],
[x^3, x^2*y, x^2*z, -x*z*y+x*z, -x*z^2, 0]]);}{%
\[
\mathit{NN} :=  \left[ 
{\begin{array}{cccccc}
x\,y & y^{2} & z\,y & 0 & 0 & z \\
0 & 0 & 0 & z\,y^{2} - z\,y & z^{2}\,y & x\,z \\
x^{2}\,z & x\,z\,y & x\,z^{2} &  - z^{2}\,y + z^{2} &  - z^{3} & 
0 \\
0 & 0 & 0 & y^{3} - 2\,y^{2} + y & z\,y^{2} - z\,y & x\,y - x \\
0 & 0 & 0 & x^{2}\,y - x^{2} - y + 1 & x^{2}\,z - z & y \\
x^{3} & x^{2}\,y & x^{2}\,z &  - x\,z\,y + x\,z &  - x\,z^{2} & 0
\end{array}}
 \right] 
\]
}
\end{maplelatex}

\end{maplegroup}
\begin{maplegroup}
The two modules $M:=\mathrm{coker}\left(R^{1\times 6}
\xrightarrow{\mathtt{M}} R^{1\times 6}\right)$ and
$N:=\mathrm{coker}\left(R^{1\times 6} \xrightarrow{\mathtt{N}}
R^{1\times 6}\right)$:

\end{maplegroup}
\begin{maplegroup}
\begin{mapleinput}
\mapleinline{active}{1d}{M:=Cokernel(MM,var);}{%
}
\end{mapleinput}

\mapleresult
\begin{maplelatex}
\mapleinline{inert}{2d}{M := [[[1, 0, 0, 0, 0, 0] = [1, 0, 0, 0, 0, 0], [0, 1, 0, 0, 0, 0] =
[0, 1, 0, 0, 0, 0], [0, 0, 1, 0, 0, 0] = [0, 0, 1, 0, 0, 0], [0, 0, 0,
1, 0, 0] = [0, 0, 0, 1, 0, 0], [0, 0, 0, 0, 1, 0] = [0, 0, 0, 0, 1,
0], [0, 0, 0, 0, 0, 1] = [0, 0, 0, 0, 0, 1]], [[x*z, z*y, z^2, 0, 0,
y], [0, 0, 0, z^2*y-z^2, z^3, x*z], [0, 0, 0, z*y^2-z*y, z^2*y, x*y],
[0, 0, 0, x*z*y-x*z, x*z^2, x^2], [x^2*z, x*z*y, x*z^2, 0, 0, x*y],
[-x*y, -y^2, -z*y, x^2*y-x^2-y+1, x^2*z-z, 0], [x^2*y-x^2, x*y^2-x*y,
x*z*y-x*z, -y^3+2*y^2-y, -z*y^2+z*y, 0], [0, 0, 0, z*y-z, z^2,
x^3-y^2]], "Presentation",
3*s^2/(1-s)+2*s/(1-s)+2/(1-s)^2+3/(1-s)^3+2*s^2/(1-s)^2+s^2+s+1/(1-s)+
s/(1-s)^2, [34, 14, 3]];}{%
\maplemultiline{
M := [[[1, \,0, \,0, \,0, \,0, \,0]=[1, \,0, \,0, \,0, \,0, \,0]
, \,[0, \,1, \,0, \,0, \,0, \,0]=[0, \,1, \,0, \,0, \,0, \,0], 
 \\
[0, \,0, \,1, \,0, \,0, \,0]=[0, \,0, \,1, \,0, \,0, \,0], \,[0, 
\,0, \,0, \,1, \,0, \,0]=[0, \,0, \,0, \,1, \,0, \,0],  \\
[0, \,0, \,0, \,0, \,1, \,0]=[0, \,0, \,0, \,0, \,1, \,0], \,[0, 
\,0, \,0, \,0, \,0, \,1]=[0, \,0, \,0, \,0, \,0, \,1]], [ \\
[x\,z, \,z\,y, \,z^{2}, \,0, \,0, \,y], \,[0, \,0, \,0, \,z^{2}\,
y - z^{2}, \,z^{3}, \,x\,z], \,[0, \,0, \,0, \,z\,y^{2} - z\,y, 
\,z^{2}\,y, \,x\,y],  \\
[0, \,0, \,0, \,x\,z\,y - x\,z, \,x\,z^{2}, \,x^{2}], \,[x^{2}\,z
, \,x\,z\,y, \,x\,z^{2}, \,0, \,0, \,x\,y],  \\
[ - x\,y, \, - y^{2}, \, - z\,y, \,x^{2}\,y - x^{2} - y + 1, \,x
^{2}\,z - z, \,0],  \\
[x^{2}\,y - x^{2}, \,x\,y^{2} - x\,y, \,x\,z\,y - x\,z, \, - y^{3
} + 2\,y^{2} - y, \, - z\,y^{2} + z\,y, \,0],  \\
[0, \,0, \,0, \,z\,y - z, \,z^{2}, \,x^{3} - y^{2}]], \,\mbox{
``Presentation''},  \\
{\displaystyle \frac {3\,s^{2}}{1 - s}}  + {\displaystyle \frac {
2\,s}{1 - s}}  + {\displaystyle \frac {2}{(1 - s)^{2}}}  + 
{\displaystyle \frac {3}{(1 - s)^{3}}}  + {\displaystyle \frac {2
\,s^{2}}{(1 - s)^{2}}}  + s^{2} + s + {\displaystyle \frac {1}{1
 - s}}  + {\displaystyle \frac {s}{(1 - s)^{2}}} ,  \\
[34, \,14, \,3]] }
}
\end{maplelatex}

\end{maplegroup}
\begin{maplegroup}
\begin{mapleinput}
\mapleinline{active}{1d}{N:=Cokernel(NN,var);}{%
}
\end{mapleinput}

\mapleresult
\begin{maplelatex}
\mapleinline{inert}{2d}{N := [[[1, 0, 0, 0, 0, 0] = [1, 0, 0, 0, 0, 0], [0, 1, 0, 0, 0, 0] =
[0, 1, 0, 0, 0, 0], [0, 0, 1, 0, 0, 0] = [0, 0, 1, 0, 0, 0], [0, 0, 0,
1, 0, 0] = [0, 0, 0, 1, 0, 0], [0, 0, 0, 0, 1, 0] = [0, 0, 0, 0, 1,
0], [0, 0, 0, 0, 0, 1] = [0, 0, 0, 0, 0, 1]], [[x*y, y^2, z*y, 0, 0,
z], [0, 0, 0, z*y^2-z*y, z^2*y, x*z], [x^2*z, x*z*y, x*z^2,
-z^2*y+z^2, -z^3, 0], [0, 0, 0, y^3-2*y^2+y, z*y^2-z*y, x*y-x], [0, 0,
0, x^2*y-x^2-y+1, x^2*z-z, y], [x^2*y, x*y^2, x*z*y, 0, 0, x*z], [x^3,
x^2*y, x^2*z, -x*z*y+x*z, -x*z^2, 0], [0, 0, 0, x*y^2*z-x*z*y,
x*y*z^2, x^2*z], [0, 0, 0, 0, 0, x^3*z-z*y^2-x*z], [0, 0, 0,
x*y^3-2*x*y^2+x*y, x*y^2*z-x*z*y, x^2*y-x^2], [0, 0, 0, 0, 0,
x^3*y-x^3-y^3-x*y+y^2+x]], "Presentation",
3*s/(1-s)+2/(1-s)^2+2*s^2+3/(1-s)^3+2*s^2/(1-s)^2+s^2/(1-s)+s^3+1/(1-s
)+s^3/(1-s)+s/(1-s)^2, [51, 17, 3]];}{%
\maplemultiline{
N := [[[1, \,0, \,0, \,0, \,0, \,0]=[1, \,0, \,0, \,0, \,0, \,0]
, \,[0, \,1, \,0, \,0, \,0, \,0]=[0, \,1, \,0, \,0, \,0, \,0], 
 \\
[0, \,0, \,1, \,0, \,0, \,0]=[0, \,0, \,1, \,0, \,0, \,0], \,[0, 
\,0, \,0, \,1, \,0, \,0]=[0, \,0, \,0, \,1, \,0, \,0],  \\
[0, \,0, \,0, \,0, \,1, \,0]=[0, \,0, \,0, \,0, \,1, \,0], \,[0, 
\,0, \,0, \,0, \,0, \,1]=[0, \,0, \,0, \,0, \,0, \,1]], [ \\
[x\,y, \,y^{2}, \,z\,y, \,0, \,0, \,z], \,[0, \,0, \,0, \,z\,y^{2
} - z\,y, \,z^{2}\,y, \,x\,z], \,[x^{2}\,z, \,x\,z\,y, \,x\,z^{2}
, \, - z^{2}\,y + z^{2}, \, - z^{3}, \,0],  \\
[0, \,0, \,0, \,y^{3} - 2\,y^{2} + y, \,z\,y^{2} - z\,y, \,x\,y
 - x], \,[0, \,0, \,0, \,x^{2}\,y - x^{2} - y + 1, \,x^{2}\,z - z
, \,y],  \\
[x^{2}\,y, \,x\,y^{2}, \,x\,z\,y, \,0, \,0, \,x\,z], \,[x^{3}, \,
x^{2}\,y, \,x^{2}\,z, \, - x\,z\,y + x\,z, \, - x\,z^{2}, \,0], 
 \\
[0, \,0, \,0, \,x\,y^{2}\,z - x\,z\,y, \,x\,y\,z^{2}, \,x^{2}\,z]
, \,[0, \,0, \,0, \,0, \,0, \,x^{3}\,z - z\,y^{2} - x\,z],  \\
[0, \,0, \,0, \,x\,y^{3} - 2\,x\,y^{2} + x\,y, \,x\,y^{2}\,z - x
\,z\,y, \,x^{2}\,y - x^{2}],  \\
[0, \,0, \,0, \,0, \,0, \,x^{3}\,y - x^{3} - y^{3} - x\,y + y^{2}
 + x]], \,\mbox{``Presentation''}, {\displaystyle \frac {3\,s}{1
 - s}}  + {\displaystyle \frac {2}{(1 - s)^{2}}}  + 2\,s^{2} \\
\mbox{} + {\displaystyle \frac {3}{(1 - s)^{3}}}  + 
{\displaystyle \frac {2\,s^{2}}{(1 - s)^{2}}}  + {\displaystyle 
\frac {s^{2}}{1 - s}}  + s^{3} + {\displaystyle \frac {1}{1 - s}
}  + {\displaystyle \frac {s^{3}}{1 - s}}  + {\displaystyle 
\frac {s}{(1 - s)^{2}}} , \,[51, \,17, \,3]] }
}
\end{maplelatex}

\end{maplegroup}
\begin{maplegroup}
The torsion submodule $K:=t(M)$ of $M$ (and $N$):

\end{maplegroup}
\begin{maplegroup}
\begin{mapleinput}
\mapleinline{active}{1d}{K:=TorsionSubmodule(M,var);}{%
}
\end{mapleinput}

\mapleresult
\begin{maplelatex}
\mapleinline{inert}{2d}{K := [[[1, 0, 0] = [0, 0, 0, 0, 0, 1], [0, 1, 0] = [0, 0, 0, y-1, z,
0], [0, 0, 1] = [x, y, z, 0, 0, 0]], [[y, 0, z], [0, -z*y, x*z], [x*z,
z^2, 0], [0, -y^2+y, x*y-x], [x*y, z*y, 0], [0, x^2-1, -y], [x^2, x*z,
0]], "Presentation", s+1/(1-s)+s/(1-s)^2+2/(1-s)^2+s/(1-s), [7, 3,
0]];}{%
\maplemultiline{
K := [[[1, \,0, \,0]=[0, \,0, \,0, \,0, \,0, \,1], \,[0, \,1, \,0
]=[0, \,0, \,0, \,y - 1, \,z, \,0], \,[0, \,0, \,1]=[x, \,y, \,z
, \,0, \,0, \,0] \\
], [[y, \,0, \,z], \,[0, \, - z\,y, \,x\,z], \,[x\,z, \,z^{2}, \,
0], \,[0, \, - y^{2} + y, \,x\,y - x], \,[x\,y, \,z\,y, \,0], \,[
0, \,x^{2} - 1, \, - y],  \\
[x^{2}, \,x\,z, \,0]], \,\mbox{``Presentation''}, \,s + 
{\displaystyle \frac {1}{1 - s}}  + {\displaystyle \frac {s}{(1
 - s)^{2}}}  + {\displaystyle \frac {2}{(1 - s)^{2}}}  + 
{\displaystyle \frac {s}{1 - s}} , \,[7, \,3, \,0]] }
}
\end{maplelatex}

\end{maplegroup}
\begin{maplegroup}
The embedding $0\to K \xrightarrow{\alpha_1} M$:

\end{maplegroup}
\begin{maplegroup}
\begin{mapleinput}
\mapleinline{active}{1d}{alpha1:=TorsionSubmoduleEmb(M,var);}{%
}
\end{mapleinput}

\mapleresult
\begin{maplelatex}
\mapleinline{inert}{2d}{alpha1 := matrix([[0, 0, 0, 0, 0, 1], [0, 0, 0, y-1, z, 0], [x, y, z,
0, 0, 0]]);}{%
\[
\alpha 1 :=  \left[ 
{\begin{array}{cccccr}
0 & 0 & 0 & 0 & 0 & 1 \\
0 & 0 & 0 & y - 1 & z & 0 \\
x & y & z & 0 & 0 & 0
\end{array}}
 \right] 
\]
}
\end{maplelatex}

\end{maplegroup}
\begin{maplegroup}
The torsion free factor $L:=M/K=M/t(M)$ of $M$ (and $N$):

\end{maplegroup}
\begin{maplegroup}
\begin{mapleinput}
\mapleinline{active}{1d}{L:=Cokernel(alpha1,M,var);}{%
}
\end{mapleinput}

\mapleresult
\begin{maplelatex}
\mapleinline{inert}{2d}{L := [[[1, 0, 0, 0, 0] = [1, 0, 0, 0, 0, 0], [0, 1, 0, 0, 0] = [0, 1,
0, 0, 0, 0], [0, 0, 1, 0, 0] = [0, 0, 1, 0, 0, 0], [0, 0, 0, 1, 0] =
[0, 0, 0, 1, 0, 0], [0, 0, 0, 0, 1] = [0, 0, 0, 0, 1, 0]], [[0, 0, 0,
y-1, z], [x, y, z, 0, 0]], "Presentation", 2/(1-s)^2+3/(1-s)^3, [5, 5,
3]];}{%
\maplemultiline{
L := [[[1, \,0, \,0, \,0, \,0]=[1, \,0, \,0, \,0, \,0, \,0], \,[0
, \,1, \,0, \,0, \,0]=[0, \,1, \,0, \,0, \,0, \,0],  \\
[0, \,0, \,1, \,0, \,0]=[0, \,0, \,1, \,0, \,0, \,0], \,[0, \,0, 
\,0, \,1, \,0]=[0, \,0, \,0, \,1, \,0, \,0],  \\
[0, \,0, \,0, \,0, \,1]=[0, \,0, \,0, \,0, \,1, \,0]], \,[[0, \,0
, \,0, \,y - 1, \,z], \,[x, \,y, \,z, \,0, \,0]], \,\mbox{
``Presentation''},  \\
{\displaystyle \frac {2}{(1 - s)^{2}}}  + {\displaystyle \frac {3
}{(1 - s)^{3}}} , \,[5, \,5, \,3]] }
}
\end{maplelatex}

\end{maplegroup}
\begin{maplegroup}
The natural epimorphism $M\xrightarrow{\alpha_2} L\to 0$:

\end{maplegroup}
\begin{maplegroup}
\begin{mapleinput}
\mapleinline{active}{1d}{alpha2:=CokernelEpi(alpha1,M,var);}{%
}
\end{mapleinput}

\mapleresult
\begin{maplelatex}
\mapleinline{inert}{2d}{alpha2 := matrix([[1, 0, 0, 0, 0], [0, 1, 0, 0, 0], [0, 0, 1, 0, 0],
[0, 0, 0, 1, 0], [0, 0, 0, 0, 1], [0, 0, 0, 0, 0]]);}{%
\[
\alpha 2 :=  \left[ 
{\begin{array}{rrrrr}
1 & 0 & 0 & 0 & 0 \\
0 & 1 & 0 & 0 & 0 \\
0 & 0 & 1 & 0 & 0 \\
0 & 0 & 0 & 1 & 0 \\
0 & 0 & 0 & 0 & 1 \\
0 & 0 & 0 & 0 & 0
\end{array}}
 \right] 
\]
}
\end{maplelatex}

\end{maplegroup}
\begin{maplegroup}
The short sequence $S: 0 \to K \xrightarrow{\alpha_1} M
\xrightarrow{\alpha_2} L \to 0$ is exact:

\end{maplegroup}
\begin{maplegroup}
\begin{mapleinput}
\mapleinline{active}{1d}{IsShortExactSeq(K,alpha1,M,alpha2,L,var);}{%
}
\end{mapleinput}

\mapleresult
\begin{maplelatex}
\mapleinline{inert}{2d}{true;}{%
\[
\mathit{true}
\]
}
\end{maplelatex}

\end{maplegroup}
\begin{maplegroup}
The functor $\Ext^1_R(-,N)$ applied to the short exact sequence $S$:

\end{maplegroup}
\begin{maplegroup}
\begin{mapleinput}
\mapleinline{active}{1d}{coseqN:=ExtOnSeqs(1,N,[alpha2,alpha1],[L,M,K],var):}{%
}
\end{mapleinput}

\end{maplegroup}
\begin{maplegroup}
This says that the map $\Ext_R^1(L,N)
\xrightarrow{\Ext_R^1(\alpha_2,N)} \Ext_R^1(M,N)$ is not injective,
i.e.\  that the connecting homomorphism
$\Hom_R(K,N)\xrightarrow{\delta^0}\Ext_R^1(L,N)$ is non-trivial:

\end{maplegroup}
\begin{maplegroup}
\begin{mapleinput}
\mapleinline{active}{1d}{IsShortExactSeq(op(MakeCoseq(coseqN)),var,"VERBOSE");}{%
}
\end{mapleinput}

\mapleresult
\begin{maplelatex}
\mapleinline{inert}{2d}{"homs" = true, "cmps" = true, "defs" = [[[1 = matrix([[0, 0, 0, 0, 0,
0], [0, 0, 0, 0, 0, 1]])], [z, y, x], "Presentation", 1, [0, 0, 0]],
true, true];}{%
\maplemultiline{
\mbox{``homs''}=\mathit{true}, \,\mbox{``cmps''}=\mathit{true}, 
\mbox{``defs''}= \\
[[[1= \left[ 
{\begin{array}{rrrrrr}
0 & 0 & 0 & 0 & 0 & 0 \\
0 & 0 & 0 & 0 & 0 & 1
\end{array}}
 \right] ], \,[z, \,y, \,x], \,\mbox{``Presentation''}, \,1, \,[0
, \,0, \,0]], \,\mathit{true}, \,\mathit{true}] }
}
\end{maplelatex}

\end{maplegroup}
\begin{maplegroup}
The functor $\Ext^1_R(-,M)$ applied to the short exact sequence $S$:

\end{maplegroup}
\begin{maplegroup}
\begin{mapleinput}
\mapleinline{active}{1d}{coseqM:=ExtOnSeqs(1,M,[alpha2,alpha1],[L,M,K],var):}{%
}
\end{mapleinput}

\end{maplegroup}
\begin{maplegroup}
Whereas the resulting sequence $\Ext_R^1(S,M)$ is again exact. In
particular, the connecting homomorphism
$\Hom_R(K,M)\xrightarrow{\delta^0}\Ext_R^1(L,M)$ is trivial, i.e.\ 
the zero map:

\end{maplegroup}
\begin{maplegroup}
\begin{mapleinput}
\mapleinline{active}{1d}{IsShortExactSeq(op(MakeCoseq(coseqM)),var,"VERBOSE");}{%
}
\end{mapleinput}

\mapleresult
\begin{maplelatex}
\mapleinline{inert}{2d}{true;}{%
\[
\mathit{true}
\]
}
\end{maplelatex}

\end{maplegroup}
%
\POP

\subsection{Example 3}\label{ExOre}

\PUSH{DerShift.tex}%
\DefineParaStyle{Maple Output}
\DefineCharStyle{2D Math}
\DefineCharStyle{2D Output}
\DefineCharStyle{LaTeX}
\begin{maplegroup}
This example demonstrates by a computation over a noncommutative {\sc
Ore} domain how one can use $\homalg$ to solve a linear system of
differential-difference equations:

\end{maplegroup}
\begin{maplegroup}
\begin{mapleinput}
\mapleinline{active}{1d}{restart:}{%
}
\end{mapleinput}

\end{maplegroup}
\begin{maplegroup}
\begin{mapleinput}
\mapleinline{active}{1d}{with(JanetOre): with(homalg):}{%
}
\end{mapleinput}

\end{maplegroup}
\begin{maplegroup}
Specify the {\tt homalg}-table of the ring package {\tt JanetOre}:

\end{maplegroup}
\begin{maplegroup}
\begin{mapleinput}
\mapleinline{active}{1d}{RPO:=`JanetOre/homalg`;}{%
}
\end{mapleinput}

\mapleresult
\begin{maplelatex}
\mapleinline{inert}{2d}{RPO := `JanetOre/homalg`;}{%
\[
\mathit{RPO} := \mathit{JanetOre/homalg}
\]
}
\end{maplelatex}

\end{maplegroup}
\begin{maplegroup}
Use the ring package {\tt JanetOre} as the default ring package:

\end{maplegroup}
\begin{maplegroup}
\begin{mapleinput}
\mapleinline{active}{1d}{`homalg/default`:=RPO;}{%
}
\end{mapleinput}

\mapleresult
\begin{maplelatex}
\mapleinline{inert}{2d}{`homalg/default` := `JanetOre/homalg`;}{%
\[
\mathit{homalg/default} := \mathit{JanetOre/homalg}
\]
}
\end{maplelatex}

\end{maplegroup}
\begin{maplegroup}
Ask homalg to check if a stably free module with a {\em finite free}
resolution\footnote{Without a finite free resolution $\homalg$ cannot
check stably freeness, cf.~\cite[Remark~28]{QR_S}.} is free,
cf.~\cite[Remark~51]{QR_S}:

\end{maplegroup}
\begin{maplegroup}
\begin{mapleinput}
\mapleinline{active}{1d}{homalg_options("check_rank_1"=true);}{%
}
\end{mapleinput}

\mapleresult
\begin{maplelatex}
\mapleinline{inert}{2d}{"JanetOre/homalg", "check_if_stably_free_rank_1_FFR_is_free" =
true;}{%
\[
\mbox{``JanetOre/homalg''}, \,\mbox{
``check\_if\_stably\_free\_rank\_1\_FFR\_is\_free''}=\mathit{true
}
\]
}
\end{maplelatex}

\end{maplegroup}
\begin{maplegroup}
Define the {\sc Ore} domain $R=\Q[t][D,t\mapsto
t,\frac{d}{dt}][\delta,t\mapsto t+1,0]$:

\end{maplegroup}
\begin{maplegroup}
\begin{mapleinput}
\mapleinline{active}{1d}{Ore:=[[t,D,delta],[],[weyl(D,t),shift(delta,t)]];}{%
}
\end{mapleinput}

\mapleresult
\begin{maplelatex}
\mapleinline{inert}{2d}{Ore := [[t, D, delta], [], [weyl(D,t), shift(delta,t)]];}{%
\[
\mathit{Ore} := [[t, \,\mathrm{D}, \,\delta ], \,[], \,[\mathrm{
weyl}(\mathrm{D}, \,t), \,\mathrm{shift}(\delta , \,t)]]
\]
}
\end{maplelatex}

\end{maplegroup}
\begin{maplegroup}
Find the general solution to the following linear system of
differential-difference equations:
\begin{eqnarray*} u \left( t \right) +u \left( t+2 \right) +\mbox {D}
\left( v \right)  \left( t+1 \right) -v \left( t+2 \right) +w \left(
t+1 \right) &=& 0, \\ t\mbox {D} \left(u \right)  \left( t+1 \right)
+u \left( t+3 \right) -2\,tu \left( t+1 \right) +2\,\mbox {D} \left( u
\right)  \left( t+1 \right) -2\,u \left( t+1 \right) \\ -v \left( t+3
\right) +2\,v \left( t+2 \right) +w \left( t+2 \right) &=& 0.
\end{eqnarray*}

\end{maplegroup}
\begin{maplegroup}
This system can be written as the following matrix operators $A$
applied to the section $\left(\begin{array}{ccc}u(t) & v(t) &
w(t)\end{array}\right)^\mathrm{tr}$. (Caution: A monomial in $t$, $D$
and $\delta$ must be read as a monomial with powers of $t$ on the left
of a monomial in the commuting $D$ and $\delta$: read in the following
both $tD$, $Dt$ as $tD$, and both $t\delta$, $\delta t$ as $t\delta$):

\end{maplegroup}
\begin{maplegroup}
\begin{mapleinput}
\mapleinline{active}{1d}{A := matrix([[1+delta^2, delta*D-delta^2, delta],
[t*delta*D+delta^3-2*t*delta+2*delta*D-2*delta, -delta^3+2*delta^2,
delta^2]]);}{%
}
\end{mapleinput}

\mapleresult
\begin{maplelatex}
\mapleinline{inert}{2d}{A := matrix([[1+delta^2, delta*D-delta^2, delta],
[t*delta*D+delta^3-2*t*delta+2*delta*D-2*delta, -delta^3+2*delta^2,
delta^2]]);}{%
\[
A :=  \left[ 
{\begin{array}{ccc}
1 + \delta ^{2} & \delta \,\mathrm{D} - \delta ^{2} & \delta  \\
t\,\delta \,\mathrm{D} + \delta ^{3} - 2\,t\,\delta  + 2\,\delta 
\,\mathrm{D} - 2\,\delta  &  - \delta ^{3} + 2\,\delta ^{2} & 
\delta ^{2}
\end{array}}
 \right] 
\]
}
\end{maplelatex}

\end{maplegroup}
\begin{maplegroup}
Studying the system is equivalent to studying the cokernel $M$ of
$A\in R^{m\times n}$ viewed as a morphism $R^{1\times
m}\xrightarrow{A} R^{1\times n}$. This is based on the following
observation: Let $\F$ be a function space, where one wants to search
for solutions of the system, then
$\Hom_R(\coker(A),\F)\cong\mathrm{Sol}_\F(A):=\{\eta\in\F^{n\times
1}\mid A\eta=0\}$, where $n$ is the number of columns of $A$.

\end{maplegroup}
\begin{maplegroup}
\begin{mapleinput}
\mapleinline{active}{1d}{M:=Cokernel(A,Ore);}{%
}
\end{mapleinput}

\mapleresult
\begin{maplelatex}
\mapleinline{inert}{2d}{M := [[[1, 0, 0] = [1, 0, 0], [0, 1, 0] = [0, 1, 0], [0, 0, 1] = [0,
0, 1]], [[1+delta^2, delta*D-delta^2, delta],
[t*delta*D+delta^3-2*t*delta+2*delta*D-2*delta, -delta^3+2*delta^2,
delta^2]], "Presentation",
3+9*s+17*s^2+s^3*(17/(1-s)+8/(1-s)^2+1/((1-s)^3)), [17, 8, 1]];}{%
\maplemultiline{
M := [[[1, \,0, \,0]=[1, \,0, \,0], \,[0, \,1, \,0]=[0, \,1, \,0]
, \,[0, \,0, \,1]=[0, \,0, \,1]],  \\
[[1 + \delta ^{2}, \,\delta \,\mathrm{D} - \delta ^{2}, \,\delta 
], \,[t\,\delta \,\mathrm{D} + \delta ^{3} - 2\,t\,\delta  + 2\,
\delta \,\mathrm{D} - 2\,\delta , \, - \delta ^{3} + 2\,\delta ^{
2}, \,\delta ^{2}]],  \\
\mbox{``Presentation''}, \,3 + 9\,s + 17\,s^{2} + s^{3}\,(
{\displaystyle \frac {17}{1 - s}}  + {\displaystyle \frac {8}{(1
 - s)^{2}}}  + {\displaystyle \frac {1}{(1 - s)^{3}}} ), \,[17, 
\,8, \,1]] }
}
\end{maplelatex}

\end{maplegroup}
\begin{maplegroup}
The torsion free factor $F$ of $M$ turns out to be free of rank $1$
(the above mentioned option \verb+"check_rank_1"+ was used):

\end{maplegroup}
\begin{maplegroup}
\begin{mapleinput}
\mapleinline{active}{1d}{F:=TorsionFreeFactor(M,Ore);}{%
}
\end{mapleinput}

\mapleresult
\begin{maplelatex}
\mapleinline{inert}{2d}{F := [[1 = [-delta, delta-D, -1]], [0], "Presentation", 1/((1-s)^3),
[0, 0, 1]];}{%
\[
F := [[1=[ - \delta , \,\delta  - \mathrm{D}, \,-1]], \,[0], \,
\mbox{``Presentation''}, \,{\displaystyle \frac {1}{(1 - s)^{3}}
} , \,[0, \,0, \,1]]
\]
}
\end{maplelatex}

\end{maplegroup}
\begin{maplegroup}
The natural epimorphsim $M\xrightarrow{\nu}F$:

\end{maplegroup}
\begin{maplegroup}
\begin{mapleinput}
\mapleinline{active}{1d}{nu:=TorsionFreeFactorEpi(M,Ore);}{%
}
\end{mapleinput}

\mapleresult
\begin{maplelatex}
\mapleinline{inert}{2d}{nu := matrix([[delta], [t], [-2+delta+t*delta-t*D-delta^2]]);}{%
\[
\nu  :=  \left[ 
{\begin{array}{c}
\delta  \\
t \\
 - 2 + \delta  + t\,\delta  - t\,\mathrm{D} - \delta ^{2}
\end{array}}
 \right] 
\]
}
\end{maplelatex}

\end{maplegroup}
\begin{maplegroup}
Since $F$ is free one can use {\tt Leftinverse} to compute a split
$M\xrightarrow{\chi}F$ (cf.~\ref{li},(1)):

\end{maplegroup}
\begin{maplegroup}
\begin{mapleinput}
\mapleinline{active}{1d}{chi:=Leftinverse(M,nu,F,Ore);}{%
}
\end{mapleinput}

\mapleresult
\begin{maplelatex}
\mapleinline{inert}{2d}{chi := matrix([[-delta, delta-D, -1]]);}{%
\[
\chi  :=  \left[ 
{\begin{array}{ccr}
 - \delta  & \delta  - \mathrm{D} & -1
\end{array}}
 \right] 
\]
}
\end{maplelatex}

\end{maplegroup}
\begin{maplegroup}
Now compute the torsion submodule $T$ of $M$:

\end{maplegroup}
\begin{maplegroup}
\begin{mapleinput}
\mapleinline{active}{1d}{T:=TorsionSubmodule(M,Ore);}{%
}
\end{mapleinput}

\mapleresult
\begin{maplelatex}
\mapleinline{inert}{2d}{T := [[[1, 0] = [t+1, -delta, 0], [0, 1] = [-2*delta,
-t*delta+t*D+delta^2+1, t]], [[1+delta^2, delta], [-2+D,
-2*delta+delta*D], [-2*delta+delta*D, 0]], "Presentation",
2+6*s+s^2*(6/(1-s)+3/(1-s)^2), [6, 3, 0]];}{%
\maplemultiline{
T := [[[1, \,0]=[t + 1, \, - \delta , \,0], \,[0, \,1]=[ - 2\,
\delta , \, - t\,\delta  + t\,\mathrm{D} + \delta ^{2} + 1, \,t]]
,  \\
[[1 + \delta ^{2}, \,\delta ], \,[ - 2 + \mathrm{D}, \, - 2\,
\delta  + \delta \,\mathrm{D}], \,[ - 2\,\delta  + \delta \,
\mathrm{D}, \,0]], \,\mbox{``Presentation''},  \\
2 + 6\,s + s^{2}\,({\displaystyle \frac {6}{1 - s}}  + 
{\displaystyle \frac {3}{(1 - s)^{2}}} ), \,[6, \,3, \,0]] }
}
\end{maplelatex}

\end{maplegroup}
\begin{maplegroup}
And the natural embedding $T\xrightarrow{\iota}M$:

\end{maplegroup}
\begin{maplegroup}
\begin{mapleinput}
\mapleinline{active}{1d}{iota:=TorsionSubmoduleEmb(M,Ore);}{%
}
\end{mapleinput}

\mapleresult
\begin{maplelatex}
\mapleinline{inert}{2d}{iota := matrix([[t+1, -delta, 0], [-2*delta, -t*delta+t*D+delta^2+1,
t]]);}{%
\[
\iota  :=  \left[ 
{\begin{array}{ccc}
t + 1 &  - \delta  & 0 \\
 - 2\,\delta  &  - t\,\delta  + t\,\mathrm{D} + \delta ^{2} + 1
 & t
\end{array}}
 \right] 
\]
}
\end{maplelatex}

\end{maplegroup}
\begin{maplegroup}
After a few tests we found an element $\mu\in T$

\end{maplegroup}
\begin{maplegroup}
\begin{mapleinput}
\mapleinline{active}{1d}{mu:=[[delta,1]];}{%
}
\end{mapleinput}

\mapleresult
\begin{maplelatex}
\mapleinline{inert}{2d}{mu := [[delta, 1]];}{%
\[
\mu  := [[\delta , \,1]]
\]
}
\end{maplelatex}

\end{maplegroup}
\begin{maplegroup}
which turned out to be a cyclic generator of $T$ (below one identifies
$\mu$ with the morphism $R^{1\times 1}\to T:\  1\mapsto \mu$):

\end{maplegroup}
\begin{maplegroup}
\begin{mapleinput}
\mapleinline{active}{1d}{IsSurjective(mu,T,Ore);}{%
}
\end{mapleinput}

\mapleresult
\begin{maplelatex}
\mapleinline{inert}{2d}{true;}{%
\[
\mathit{true}
\]
}
\end{maplelatex}

\end{maplegroup}
\begin{maplegroup}
$H$ is nothing but $T$ rewritten on this cyclic generator. It further
turns out that the cyclic generator $\mu$ of the torsion submodule $T$
satisfies a single simple relation $\delta^2(D-2)\mu=0$:

\end{maplegroup}
\begin{maplegroup}
\begin{mapleinput}
\mapleinline{active}{1d}{H:=Image(mu,T,Ore);}{%
}
\end{mapleinput}

\mapleresult
\begin{maplelatex}
\mapleinline{inert}{2d}{H := [[1 = [t*delta, t*D-t*delta+1, t]], [-2*delta^2+delta^2*D],
"Presentation", 1+3*s+6*s^2+s^3*(6/(1-s)+3/(1-s)^2), [6, 3, 0]];}{%
\maplemultiline{
H := [[1=[t\,\delta , \,t\,\mathrm{D} - t\,\delta  + 1, \,t]], \,
[ - 2\,\delta ^{2} + \delta ^{2}\,\mathrm{D}], \,\mbox{
``Presentation''},  \\
1 + 3\,s + 6\,s^{2} + s^{3}\,({\displaystyle \frac {6}{1 - s}} 
 + {\displaystyle \frac {3}{(1 - s)^{2}}} ), \,[6, \,3, \,0]] }
}
\end{maplelatex}

\end{maplegroup}
\begin{maplegroup}
$\varepsilon$ is nothing but $\iota$ rewritten on the cyclic generator
$\mu$:

\end{maplegroup}
\begin{maplegroup}
\begin{mapleinput}
\mapleinline{active}{1d}{epsilon:=ComposeMaps(mu,iota,M,Ore);}{%
}
\end{mapleinput}

\mapleresult
\begin{maplelatex}
\mapleinline{inert}{2d}{epsilon := matrix([[t*delta, t*D-t*delta+1, t]]);}{%
\[
\varepsilon  :=  \left[ 
{\begin{array}{ccc}
t\,\delta  & t\,\mathrm{D} - t\,\delta  + 1 & t
\end{array}}
 \right] 
\]
}
\end{maplelatex}

\end{maplegroup}
\begin{maplegroup}
\begin{mapleinput}
\mapleinline{active}{1d}{IsHom(H,epsilon,M,Ore);}{%
}
\end{mapleinput}

\mapleresult
\begin{maplelatex}
\mapleinline{inert}{2d}{true;}{%
\[
\mathit{true}
\]
}
\end{maplelatex}

\end{maplegroup}
\begin{maplegroup}
\begin{mapleinput}
\mapleinline{active}{1d}{IsInjective(H,epsilon,M,Ore);}{%
}
\end{mapleinput}

\mapleresult
\begin{maplelatex}
\mapleinline{inert}{2d}{true;}{%
\[
\mathit{true}
\]
}
\end{maplelatex}

\end{maplegroup}
\begin{maplegroup}
Now construct the isomorphism $\alpha$ from the direct sum $S:=F\oplus
H$ onto $M$:

\end{maplegroup}
\begin{maplegroup}
\begin{mapleinput}
\mapleinline{active}{1d}{alpha:=RPO[matrix](RPO[UnionOfRows](chi,epsilon));}{%
}
\end{mapleinput}

\mapleresult
\begin{maplelatex}
\mapleinline{inert}{2d}{alpha := matrix([[-delta, delta-D, -1], [t*delta, t*D-t*delta+1,
t]]);}{%
\[
\alpha  :=  \left[ 
{\begin{array}{ccc}
 - \delta  & \delta  - \mathrm{D} & -1 \\
t\,\delta  & t\,\mathrm{D} - t\,\delta  + 1 & t
\end{array}}
 \right] 
\]
}
\end{maplelatex}

\end{maplegroup}
\begin{maplegroup}
\begin{mapleinput}
\mapleinline{active}{1d}{S:=DirectSum(F,H,Ore);}{%
}
\end{mapleinput}

\mapleresult
\begin{maplelatex}
\mapleinline{inert}{2d}{S := [[[1, 0] = [1, 0], [0, 1] = [0, 1]], [[0,
-2*delta^2+delta^2*D]], "Presentation",
2+6*s+12*s^2+s^3*(12/(1-s)+6/(1-s)^2+1/((1-s)^3)), [12, 6, 1]];}{%
\maplemultiline{
S := [[[1, \,0]=[1, \,0], \,[0, \,1]=[0, \,1]], \,[[0, \, - 2\,
\delta ^{2} + \delta ^{2}\,\mathrm{D}]], \,\mbox{``Presentation''
},  \\
2 + 6\,s + 12\,s^{2} + s^{3}\,({\displaystyle \frac {12}{1 - s}} 
 + {\displaystyle \frac {6}{(1 - s)^{2}}}  + {\displaystyle 
\frac {1}{(1 - s)^{3}}} ), \,[12, \,6, \,1]] }
}
\end{maplelatex}

\end{maplegroup}
\begin{maplegroup}
\begin{mapleinput}
\mapleinline{active}{1d}{IsHom(S,alpha,M,Ore);}{%
}
\end{mapleinput}

\mapleresult
\begin{maplelatex}
\mapleinline{inert}{2d}{true;}{%
\[
\mathit{true}
\]
}
\end{maplelatex}

\end{maplegroup}
\begin{maplegroup}
\begin{mapleinput}
\mapleinline{active}{1d}{IsBijective(S,alpha,M,Ore);}{%
}
\end{mapleinput}

\mapleresult
\begin{maplelatex}
\mapleinline{inert}{2d}{true;}{%
\[
\mathit{true}
\]
}
\end{maplelatex}

\end{maplegroup}
\begin{maplegroup}
$N$ is now nothing but $M$ rewritten on the images of the two
generators of $S$. So $M$ (or $N$) is a module generated by two
generators, one free and one subject to a single simple relation:

\end{maplegroup}
\begin{maplegroup}
\begin{mapleinput}
\mapleinline{active}{1d}{N:=Image(alpha,M,Ore);}{%
}
\end{mapleinput}

\mapleresult
\begin{maplelatex}
\mapleinline{inert}{2d}{N := [[[1, 0] = [-delta, delta-D, -1], [0, 1] = [t*delta,
t*D-t*delta+1, t]], [[0, -2*delta^2+delta^2*D]], "Presentation",
2+6*s+12*s^2+s^3*(12/(1-s)+6/(1-s)^2+1/((1-s)^3)), [12, 6, 1]];}{%
\maplemultiline{
N := [[[1, \,0]=[ - \delta , \,\delta  - \mathrm{D}, \,-1], \,[0
, \,1]=[t\,\delta , \,t\,\mathrm{D} - t\,\delta  + 1, \,t]], \,[[
0, \, - 2\,\delta ^{2} + \delta ^{2}\,\mathrm{D}]],  \\
\mbox{``Presentation''}, \,2 + 6\,s + 12\,s^{2} + s^{3}\,(
{\displaystyle \frac {12}{1 - s}}  + {\displaystyle \frac {6}{(1
 - s)^{2}}}  + {\displaystyle \frac {1}{(1 - s)^{3}}} ), \,[12, 
\,6, \,1]] }
}
\end{maplelatex}

\end{maplegroup}
\begin{maplegroup}
Using sections, the above presentation reads: The original system on
the three unkown functions $u(t)$, $v(t)$ and $w(t)$ is equivalent to
a system on two unknown functions $f(t)$ and $g(t)$, where only $g(t)$
satisfies the single simple equation $D(g)(t+1)-g(t+1)=0$. It is also
given how to express $f$ and $g$ in terms of $u$, $v$ and $w$.

\end{maplegroup}
\begin{maplegroup}
\begin{mapleinput}
\mapleinline{active}{1d}{PresentationOnSections(N,Ore[1],Ore[3],[u,v,w],[f,g]);}{%
}
\end{mapleinput}

\mapleresult
\begin{maplelatex}
\mapleinline{inert}{2d}{[[f(t) = -u(t+1)+v(t+1)-D(v)(t)-w(t), g(t) =
t*u(t+1)+t*D(v)(t)-t*v(t+1)+v(t)+t*w(t)], [-2*g(t+2)+D(g)(t+2)],
"Presentation", 2+6*s+12*s^2+s^3*(12/(1-s)+6/(1-s)^2+1/((1-s)^3)),
[12, 6, 1]];}{%
\maplemultiline{
[[\mathrm{f}(t)= - \mathrm{u}(t + 1) + \mathrm{v}(t + 1) - 
\mathrm{D}(v)(t) - \mathrm{w}(t),  \\
\mathrm{g}(t)=t\,\mathrm{u}(t + 1) + t\,\mathrm{D}(v)(t) - t\,
\mathrm{v}(t + 1) + \mathrm{v}(t) + t\,\mathrm{w}(t)], \,[ - 2\,
\mathrm{g}(t + 2) + \mathrm{D}(g)(t + 2)],  \\
\mbox{``Presentation''}, \,2 + 6\,s + 12\,s^{2} + s^{3}\,(
{\displaystyle \frac {12}{1 - s}}  + {\displaystyle \frac {6}{(1
 - s)^{2}}}  + {\displaystyle \frac {1}{(1 - s)^{3}}} ), \,[12, 
\,6, \,1]] }
}
\end{maplelatex}

\end{maplegroup}
\begin{maplegroup}
$f(t)$ is therefore a free function and $g(t)=C\exp(t)$, where $C$ is
an arbitrary constant:

\end{maplegroup}
\begin{maplegroup}
\begin{mapleinput}
\mapleinline{active}{1d}{sol:=[f=unapply(f(t),t),g=unapply(C*exp(2*t),t)];}{%
}
\end{mapleinput}

\mapleresult
\begin{maplelatex}
\mapleinline{inert}{2d}{sol := [f = (t -> f(t)), g = (t -> C*exp(2*t))];}{%
\[
\mathit{sol} := [f=(t\rightarrow \mathrm{f}(t)), \,g=(t
\rightarrow C\,e^{(2\,t)})]
\]
}
\end{maplelatex}

\end{maplegroup}
\begin{maplegroup}
In order to solve the original system one needs to express $u$, $v$
and $w$ in terms of $f$ and $g$. To this end use {\tt Leftinverse} to
compute $\eta:=\alpha^{-1}$ (cf.~\ref{li},(2)):

\end{maplegroup}
\begin{maplegroup}
\begin{mapleinput}
\mapleinline{active}{1d}{eta:=Leftinverse(S,alpha,M,Ore);}{%
}
\end{mapleinput}

\mapleresult
\begin{maplelatex}
\mapleinline{inert}{2d}{eta := matrix([[delta, 0], [t, 1], [-2+delta+t*delta-t*D-delta^2,
delta-D]]);}{%
\[
\eta  :=  \left[ 
{\begin{array}{cc}
\delta  & 0 \\
t & 1 \\
 - 2 + \delta  + t\,\delta  - t\,\mathrm{D} - \delta ^{2} & 
\delta  - \mathrm{D}
\end{array}}
 \right] 
\]
}
\end{maplelatex}

\end{maplegroup}
\begin{maplegroup}
Again, $L$ is nothing but $S$ ($\cong M$) rewritten on the images
(under $\eta$) of the three original generators of $M$:

\end{maplegroup}
\begin{maplegroup}
\begin{mapleinput}
\mapleinline{active}{1d}{L:=Image(eta,S,Ore,"USE_IMAGE_OF_GENERATORS");}{%
}
\end{mapleinput}

\mapleresult
\begin{maplelatex}
\mapleinline{inert}{2d}{L := [[[1, 0, 0] = [delta, 0], [0, 1, 0] = [t, 1], [0, 0, 1] =
[-2+delta+t*delta-t*D-delta^2, delta-D]], [[1+delta^2,
delta*D-delta^2, delta],
[t*delta*D+delta^3-2*t*delta+2*delta*D-2*delta, -delta^3+2*delta^2,
delta^2]], "Presentation",
3+9*s+17*s^2+s^3*(17/(1-s)+8/(1-s)^2+1/((1-s)^3)), [17, 8, 1]];}{%
\maplemultiline{
L := [[[1, \,0, \,0]=[\delta , \,0], \,[0, \,1, \,0]=[t, \,1], \,
[0, \,0, \,1]=[ - 2 + \delta  + t\,\delta  - t\,\mathrm{D} - 
\delta ^{2}, \,\delta  - \mathrm{D}]],  \\
[[1 + \delta ^{2}, \,\delta \,\mathrm{D} - \delta ^{2}, \,\delta 
], \,[t\,\delta \,\mathrm{D} + \delta ^{3} - 2\,t\,\delta  + 2\,
\delta \,\mathrm{D} - 2\,\delta , \, - \delta ^{3} + 2\,\delta ^{
2}, \,\delta ^{2}]],  \\
\mbox{``Presentation''}, \,3 + 9\,s + 17\,s^{2} + s^{3}\,(
{\displaystyle \frac {17}{1 - s}}  + {\displaystyle \frac {8}{(1
 - s)^{2}}}  + {\displaystyle \frac {1}{(1 - s)^{3}}} ), \,[17, 
\,8, \,1]] }
}
\end{maplelatex}

\end{maplegroup}
\begin{maplegroup}
And again, using the {\tt JanetOre} procedure {\tt
PresentationOnSections}, the above presentation is rewritten:

\end{maplegroup}
\begin{maplegroup}
\begin{mapleinput}
\mapleinline{active}{1d}{P:=PresentationOnSections(L,Ore[1],Ore[3],[f,g],[u,v,w]);}{%
}
\end{mapleinput}

\mapleresult
\begin{maplelatex}
\mapleinline{inert}{2d}{P := [[u(t) = f(t+1), v(t) = t*f(t)+g(t), w(t) =
-2*f(t)+f(t+1)+t*f(t+1)-t*D(f)(t)-f(t+2)+g(t+1)-D(g)(t)],
[u(t)+u(t+2)+D(v)(t+1)-v(t+2)+w(t+1),
t*D(u)(t+1)+u(t+3)-2*t*u(t+1)+2*D(u)(t+1)-2*u(t+1)-v(t+3)+2*v(t+2)+w(t
+2)], "Presentation",
3+9*s+17*s^2+s^3*(17/(1-s)+8/(1-s)^2+1/((1-s)^3)), [17, 8, 1]];}{%
\maplemultiline{
P := [[\mathrm{u}(t)=\mathrm{f}(t + 1), \,\mathrm{v}(t)=t\,
\mathrm{f}(t) + \mathrm{g}(t),  \\
\mathrm{w}(t)= - 2\,\mathrm{f}(t) + \mathrm{f}(t + 1) + t\,
\mathrm{f}(t + 1) - t\,\mathrm{D}(f)(t) - \mathrm{f}(t + 2) + 
\mathrm{g}(t + 1) - \mathrm{D}(g)(t)], [ \\
\mathrm{u}(t) + \mathrm{u}(t + 2) + \mathrm{D}(v)(t + 1) - 
\mathrm{v}(t + 2) + \mathrm{w}(t + 1), t\,\mathrm{D}(u)(t + 1) + 
\mathrm{u}(t + 3) \\
\mbox{} - 2\,t\,\mathrm{u}(t + 1) + 2\,\mathrm{D}(u)(t + 1) - 2\,
\mathrm{u}(t + 1) - \mathrm{v}(t + 3) + 2\,\mathrm{v}(t + 2) + 
\mathrm{w}(t + 2)],  \\
\mbox{``Presentation''}, \,3 + 9\,s + 17\,s^{2} + s^{3}\,(
{\displaystyle \frac {17}{1 - s}}  + {\displaystyle \frac {8}{(1
 - s)^{2}}}  + {\displaystyle \frac {1}{(1 - s)^{3}}} ), \,[17, 
\,8, \,1]] }
}
\end{maplelatex}

\end{maplegroup}
\begin{maplegroup}
The most general $C^\infty(\R)$-solution is given by:

\end{maplegroup}
\begin{maplegroup}
\begin{mapleinput}
\mapleinline{active}{1d}{Sol:=eval(RPO[matrix](GeneratorsOfPresentation(P)),sol);}{%
}
\end{mapleinput}

\mapleresult
\begin{maplelatex}
\mapleinline{inert}{2d}{Sol := matrix([[f(t+1)], [t*f(t)+C*exp(2*t)],
[-2*f(t)+f(t+1)+t*f(t+1)-t*D(f)(t)-f(t+2)+C*exp(2*t+2)-2*C*exp(2*t)]])
;}{%
\[
\mathit{Sol} :=  \left[ 
{\begin{array}{c}
\mathrm{f}(t + 1) \\
t\,\mathrm{f}(t) + C\,e^{(2\,t)} \\
 - 2\,\mathrm{f}(t) + \mathrm{f}(t + 1) + t\,\mathrm{f}(t + 1) - 
t\,\mathrm{D}(f)(t) - \mathrm{f}(t + 2) + C\,e^{(2\,t + 2)} - 2\,
C\,e^{(2\,t)}
\end{array}}
 \right] 
\]
}
\end{maplelatex}

\end{maplegroup}
\begin{maplegroup}
Besides, one has the following uniqueness property: If $\widetilde{f}$
and $\widetilde{C}$ lead to the same solution then $\widetilde{f}=f$
and $\widetilde{C}=C$.

\end{maplegroup}
\begin{maplegroup}
Test the general solution with the {\tt JanetOre} procedure {\tt
JApplyMatrix}:

\end{maplegroup}
\begin{maplegroup}
\begin{mapleinput}
\mapleinline{active}{1d}{JApplyMatrix(A,Sol,Ore[1],Ore[3]);}{%
}
\end{mapleinput}

\mapleresult
\begin{maplelatex}
\mapleinline{inert}{2d}{matrix([[0], [0]]);}{%
\[
 \left[ 
{\begin{array}{r}
0 \\
0
\end{array}}
 \right] 
\]
}
\end{maplelatex}

\end{maplegroup}
%
\POP

%
\POP

\PUSH{homalg.bbl}%
\newcommand{\etalchar}[1]{$^{#1}$}
\providecommand{\bysame}{\leavevmode\hbox to3em{\hrulefill}\thinspace}
\providecommand{\MR}{\relax\ifhmode\unskip\space\fi MR }

\providecommand{\href}[2]{#2}

\POP



\edoc
